\documentclass[11pt]{article}
\usepackage{amssymb}
\usepackage{times}

\title{Tests \`a la Hurewicz dans le plan.\indent}
\author{Dominique LECOMTE}
\date{\it ~Fund. Math.~\rm 156 (1998), 131-165}

\newcommand{\Ana}{{\it\Sigma}^{1}_{1}}

\newcommand{\Borel}{{\it\Delta}^{1}_{1}}

\newcommand{\boraone}{{\bf\Sigma}^{0}_{1}}
\newcommand{\boratwo}{{\bf\Sigma}^{0}_{2}}
\newcommand{\borathree}{{\bf\Sigma}^{0}_{3}}

\newcommand{\boraxi}{{\bf\Sigma}^{0}_{\xi}}

\newcommand{\bortwo}{{\bf\Delta}^{0}_{2}}
\newcommand{\borone}{{\bf\Delta}^{0}_{1}}

\newcommand{\bormone}{{\bf\Pi}^{0}_{1}}
\newcommand{\bormtwo}{{\bf\Pi}^{0}_{2}}
\newcommand{\bormthree}{{\bf\Pi}^{0}_{3}}

\newcommand{\boraom}{{\bf\Sigma}^{0}_{\omega}}
\newcommand{\bormom}{{\bf\Pi}^{0}_{\omega}}

\newtheorem{thm} {Th\'eor\`eme} [section]
\newtheorem{defi} [thm] {D\'efinition}
\newtheorem{defis} [thm] {D\'efinitions}

\newtheorem{lem} [thm] {Lemme}

\begin{document}

\maketitle

\noindent {\footnotesize {\bf R\'esum\'e.} Nous donnons, pour une certaine cat\'egorie de bor\'eliens d'un produit de deux espaces polonais, comprenant les bor\'eliens \`a coupes d\'enombrables, une 
caract\'erisation du type ``test d'Hurewicz" de ceux ne pouvant pas \^etre rendus  
diff\'erence transfinie d'ouverts par changement des deux topologies polonaises.}

\section{$\!\!\!\!\!\!$ Introduction.}\indent

 Ces travaux se situent dans le cadre de la th\'eorie descriptive des ensembles. Je 
renvoie le lecteur \`a [Ku] (resp. [Mo]) pour les notions de base de th\'eorie descriptive 
classique (resp. effective). Rappelons que dans le cas des 
espaces polonais de dimension 0, la hi\'erarchie de Baire des bor\'eliens est construite 
en alternant les op\'erations de r\'eunion d\'enombrable et de passage au compl\'ementaire, en 
partant des ouverts-ferm\'es, ce de mani\`ere transfinie. On a alors la hi\'erarchie suivante :
$$\begin{array}{ll} 
& \boraone=\mbox{ouverts}~~~\boratwo=F_{\sigma} ~~~...~~~\boraom ~~~...\cr 
& \bormone=\mbox{ferm\'es}~~~\bormtwo=G_{\delta} ~~~...~~~\bormom ~~~...
\end{array}$$
On s'int\'eresse ici \`a une hi\'erarchie analogue \`a celle de Baire, sauf qu'au lieu de 
partir des ouverts-ferm\'es d'un espace polonais de dimension 0, on part des produits de 
deux bor\'eliens, chacun d'entre eux \'etant inclus dans un espace polonais. L'analogie 
devient plus claire quand on sait qu'\'etant donn\'es un espace polonais $X$ et un 
bor\'elien $A$ de $X$, on peut trouver une topologie polonaise plus fine que la 
topologie initiale sur $X$ (topologie ayant donc les m\^emes bor\'eliens), de dimension 0, 
et qui rende $A$ ouvert-ferm\'e. Pour notre probl\`eme, le fait de travailler dans les 
espaces de dimension 0 n'est donc pas une restriction r\'eelle. La d\'efinition qui suit 
appara\^\i t alors naturelle :\bigskip

\noindent\bf D\'efinition.\it\ Soient $X$ et $Y$ des espaces polonais, et $A$ un bor\'elien de 
$X \times Y$. Si $\Gamma$ est une classe de Baire, on dira que 
$A$ est $potentiellement\ dans\ \Gamma$ $($ce qu'on notera $A \in 
\mbox{pot}(\Gamma))$ s'il existe des topologies polonaises de dimension 0, 
$\sigma$ $($sur $X)$ et $\tau$ $($sur $Y)$, plus fines que les topologies 
initiales, telles que $A$, consid\'er\'e comme partie de $(X, {\sigma}) \times (Y, \tau)$, soit dans $\Gamma$.\rm\bigskip

 La motivation pour l'\'etude de ces classes de Baire potentielles trouve son 
origine dans l'\'etude des relations d'\'equivalence bor\'eliennes, et plus pr\'ecis\'ement 
dans l'\'etude du pr\'e-ordre suivant sur la collection des relations d'\'equivalence 
bor\'eliennes d\'efinies sur un espace polonais : 
$$E\leq F~~\Leftrightarrow ~~\exists~f~\mbox{bor\'elienne}~~E=(f\times f)^{-1}(F).$$
A l'aide de la notion de classe de Baire potentielle, A. Louveau montre dans 
[Lo3] que la collection des relations d'\'equivalence $\boraxi$ n'est pas 
co-finale, et il en d\'eduit qu'il n'existe pas de relation maximum pour $\leq$.

\vfill\eject

 Pour d\'eterminer la complexit\'e exacte d'un bor\'elien, on est amen\'e \`a montrer qu'il n'est 
pas d'une classe de Baire donn\'ee - ce qui est g\'en\'eralement beaucoup plus difficile que 
de montrer qu'il est d'une autre classe de Baire. Le th\'eor\`eme d'Hurewicz, rappel\'e 
ci-dessous, donne une condition n\'ecessaire et suffisante pour la classe des $G_{\delta}$ 
(cf [SR]).\bigskip

\noindent\bf Th\'eor\`eme.\it\ Soient $X$ un espace polonais, et $A$ un bor\'elien de $X$. Alors les 
conditions suivantes sont \'equivalentes :\smallskip

\noindent (a) Le bor\'elien $A$ n'est pas $\bormtwo$.\smallskip

\noindent (b) Il existe $u:2^\omega\rightarrow X$ injective continue telle que 
$u^{-1}(A)=\{\alpha\in2^\omega~/~\exists~n~~\forall~m\geq n~~\alpha (m)=0\}.$\rm\bigskip

 Ce th\'eor\`eme a \'et\'e g\'en\'eralis\'e aux autres classes de Baire par A. Louveau et J. Saint Raymond (cf [Lo-SR]). On cherche \`a \'etablir des r\'esultats analogues au th\'eor\`eme d'Hurewicz pour les classes 
de Baire potentielles. Dans la premi\`ere partie, nous nous 
int\'eresseront \`a la caract\'erisation des ensembles potentiellement ferm\'es ; nous 
d\'emontrons le\bigskip

\noindent\bf Th\'eor\`eme.\it\ Il existe un bor\'elien $A_1$ de $2^\omega\times 2^\omega$ tel que pour 
tous espaces polonais $X$ et $Y$, et pour tout bor\'elien $A$ de $X\times Y$ qui est 
$\mbox{pot}(\borathree)$ et $\mbox{pot}(\bormthree)$, on a l'\'equivalence entre les conditions suivantes :\smallskip

\noindent (a) Le bor\'elien $A$ n'est pas $\mbox{pot}(\bormone)$.\smallskip

\noindent (b) Il existe des fonctions continues $u : 2^\omega \rightarrow X$ et 
$v : 2^\omega \rightarrow Y$ telles que $\overline{A_1} \cap (u\times v)^{-1}(A) = 
A_1$.\rm\bigskip

 Rappelons que les bor\'eliens \`a coupes verticales (ou horizontales) d\'enombrables sont 
$\mbox{pot}(\boratwo)$ (cf [Lo1]), donc v\'erifient l'hypoth\`ese de ce th\'eor\`eme. Il est \`a 
noter qu'on ne peut pas esp\'erer une r\'eduction sur tout le produit, c'est \`a dire qu'on 
ne peut pas avoir $(u\times v)^{-1}(A) = A_1$ dans la condition (b) (cf [Le1], 
Cor. 4.14.(b)). Nous montrerons \'egalement l'impossibilit\'e d'avoir l'injectivit\'e des 
fonctions $u$ et $v$ de r\'eduction, ce qui constitue une autre diff\'erence avec le 
th\'eor\`eme d'Hurewicz.\bigskip

 Dans la seconde partie, nous \'etendrons ce th\'eor\`eme \`a d'autres 
classes, qui s'introduisent naturellement \`a partir du th\'eor\`eme d'Hurewicz. En effet, ce 
th\'eor\`eme montre entre autres l'int\'er\^et des r\'eductions par des fonctions continues pour 
la comparaison de la complexit\'e des bor\'eliens. La d\'efinition suivante appara\^\i t alors 
naturelle :\bigskip

\noindent\bf D\'efinition.\it\ Soit $\Gamma$ une classe de parties d'espaces 
polonais de dimension 0. On dit que $\Gamma$ est une $classe\ de\ Wadge$ s'il 
existe un espace polonais $P_0$ de dimension 0, et un bor\'elien $A_0$ de 
$P_0$ tels que pour tout espace polonais $P$ de dimension 0 et pour toute 
partie $A$ de $P$, $A$ est dans $\Gamma$ si et seulement s'il existe une 
fonction continue $f$ de $P$ dans $P_0$ telle que $A = f^{-1}(A_0)$.\rm\bigskip

 On peut d\'emontrer que la hi\'erarchie de Wadge affine celle de Baire. L'utilit\'e de 
consid\'erer les espaces de dimension 0 appara\^\i t ici : il faut assurer l'existence de 
suffisamment de fonctions continues. En effet, les seules fonctions continues de $\mathbb{R}$ dans 
$2^\omega$ sont les fonctions constantes ! Il y a eu des travaux, notamment de A. 
Louveau et J. Saint Raymond (cf [Lo2]), pour d\'ecrire la hi\'erarchie de Wadge en 
termes d'op\'erations ensemblistes, comme dans la hi\'erarchie de Baire. Ceci am\`ene \`a 
consid\'erer de nouveaux ensembles. 

\vfill\eject

 Par exemple, si $\xi$ est un ordinal d\'enombrable 
et $(A_\eta)_{\eta <\xi}$ une suite croissante d'ouverts d'un espace polonais $X$, on note 
$$D((A_\eta)_{\eta <\xi}) := \left\{x\in X~/~\exists~\eta <\xi~~x\in A_\eta\setminus (
\bigcup_{\theta <\eta} A_\theta )~\mbox{et}~\eta~\mbox{n'a~pas~la~m\^eme~parit\'e~que}~\xi\right\}.$$
On note $D_\xi(\boraone)$ la classe des ensembles de la forme 
$D((A_\eta)_{\eta <\xi})$. On peut montrer que les seules classes 
de Wadge non stables par passage au compl\'ementaire contenues dans 
$\bortwo=\boratwo\cap\bormtwo$ sont les 
$D_\xi(\boraone)$ et $\check D_\xi(\boraone)$. On obtient alors la hi\'erarchie 
suivante :
$$\begin{array}{ll}  
D_0 (\boraone)=\{\emptyset\} 
&~~~D_1 (\boraone)=\mbox{ouverts} ~~~...~~~D_\omega (\boraone)~~~...~~~\boratwo \cr  
\check D_0 (\boraone) 
&~~~\check D_1 (\boraone)=\mbox{ferm\'es } ~~~...~~~\check D_\omega (\boraone)~~~...~~~\bormtwo 
\end{array}$$ 
On peut d\'efinir sans probl\`eme les ensembles potentiellement $\Gamma$, o\`u $\Gamma$ 
est une classe de Wadge, en utilisant la m\^eme d\'efinition que pr\'ec\'edemment. Nous d\'emontrons le\bigskip

\noindent\bf Th\'eor\`eme.\it\ Soit $\xi$ un ordinal d\'enombrable.\smallskip

\noindent (1) Si $\xi$ est pair, il existe un bor\'elien $A_\xi$ de $2^\omega\times 2^\omega$ tel que 
pour tous espaces polonais $X$ et $Y$, et pour tout bor\'elien $A$ de $X\times Y$ qui est 
$\mbox{pot}(\borathree)$ et $\mbox{pot}(\bormthree)$, on a l'\'equivalence entre les conditions suivantes :\smallskip

\noindent (a) Le bor\'elien $A$ n'est pas $\mbox{pot}(D_\xi (\boraone))$.\smallskip

\noindent (b) Il existe des fonctions continues $u : 2^\omega \rightarrow X$ et 
$v : 2^\omega \rightarrow Y$ telles que $\overline{A_\xi} \cap (u\times v)^{-1}(A) = A_\xi$.\medskip

\noindent (2) Si $\xi$ est impair, il existe un bor\'elien $A_\xi$ de $2^\omega\times 2^\omega$ tel que 
pour tous espaces polonais $X$ et $Y$, et pour tout bor\'elien $A$ de $X\times Y$ qui est 
$\mbox{pot}(\borathree)$ et $\mbox{pot}(\bormthree)$, on a l'\'equivalence entre les conditions suivantes :\smallskip

\noindent (a) Le bor\'elien $A$ n'est pas $\mbox{pot}(\check D_\xi (\boraone))$.\smallskip

\noindent (b) Il existe des fonctions continues $u : 2^\omega \rightarrow X$ et 
$v : 2^\omega \rightarrow Y$ telles que $\overline{A_\xi} \cap (u\times v)^{-1}(A) = A_\xi$.\rm\bigskip

\noindent\bf Questions.\rm~(a) Un premier probl\`eme ouvert est de savoir si on peut 
supprimer l'hypoth\`ese ``$A$ est $\mbox{pot}(\borathree)$ et $\mbox{pot}(\bormthree)$" dans ce th\'eor\`eme.\bigskip

\noindent (b) Un deuxi\`eme probl\`eme ouvert est le suivant. Comme nous l'avons mentionn\'e 
avant, les seules classes de Wadge non stables par passage au compl\'ementaire 
contenues dans $\bortwo$ sont les $D_\xi(\boraone)$ et $\check D_\xi(\boraone)$. Les 
bor\'eliens \`a coupes d\'enombrables \'etant $\mbox{pot}(\boratwo)$, la caract\'erisation de ces 
bor\'eliens en termes de ``tests \`a la Hurewicz" est donc compl\`ete, \`a l'exception de 
ceux qui ne sont pas $\mbox{pot}(\bormtwo)$. La question est donc de savoir si la conjecture 
suivante est vraie :\bigskip

\noindent\bf Conjecture.\it\ Il existe un bor\'elien $B$ de $2^\omega\times2^\omega$, tel que pour 
tous espaces polonais $X$ et $Y$, et pour tout bor\'elien $A$ de $X\times Y$ \`a coupes 
d\'enombrables, on a l'\'equivalence entre les conditions suivantes :\smallskip

\noindent (a) Le bor\'elien $A$ n'est pas $\mbox{pot}(\bormtwo)$.\smallskip

\noindent (b) Il existe des fonctions continues $u : 2^\omega \rightarrow X$ et 
$v : 2^\omega \rightarrow Y$ telles que $\overline{B} \cap (u\times v)^{-1}(A) = B$.\rm

\section{$\!\!\!\!\!\!$ Un test pour les ensembles non potentiellement ferm\'es.}

\noindent\bf (A) La construction de base.\rm\bigskip

 Nous montrons un premier r\'esultat qui n'est pas tout \`a fait la caract\'erisation des 
ensembles non potentiellement ferm\'es annonc\'ee dans l'introduction. Sa d\'emonstration 
est plus importante que l'\'enonc\'e lui-m\^eme, et fournit une construction qui sera affin\'ee 
plus tard, de 3 mani\`eres diff\'erentes :\bigskip

\noindent - Pour \'etablir la caract\'erisation des ensembles non potentiellement ferm\'es.\smallskip
 
\noindent - Pour montrer l'impossibilit\'e de l'injectivit\'e de la r\'eduction.\smallskip

\noindent - Pour \'etablir la caract\'erisation des ensembles non potentiellement diff\'erence 
transfinie d'ouverts.\bigskip

 Pour \'enoncer et \'etablir ce r\'esultat, il nous faut du vocabulaire :
  
\begin{defi} Soient $(G_n)$ une suite de ferm\'es et $G$ un ferm\'e d'un espace 
topologique $X$. On dit que $(G_n)\ converge\ vers\ G$ si 
$G=\overline{\bigcup_{n\in\omega} G_n}\setminus (\bigcup_{n\in\omega} G_n)$.\end{defi}

 L'id\'ee est la suivante : comment tester si une partie $A$ d'un espace m\'etrique est 
ferm\'ee ? Une r\'eponse est que $A$ n'est pas ferm\'e si et seulement s'il existe une suite 
$(x_n)$ d'\'el\'ements de $A$ convergeant vers un point $x$ hors de $A$. On a alors, avec la 
d\'efinition pr\'ec\'edente, que $(\{x_n\})$ converge vers $\{ x\}$. On ne peut pas prendre ce 
test pour caract\'eriser les ensembles non potentiellement ferm\'es, puisque le singleton  
$\{ x\}$ peut \^etre rendu ouvert-ferm\'e. Cependant, on peut remarquer que si $X$ est 
un espace polonais et $\tau$ une topologie polonaise plus fine sur $X$, il existe un 
$G_{\delta}$ dense de $X$ sur lequel les deux topologies co\"\i ncident. L'id\'ee est donc de 
remplacer les singletons par des ensembles rencontrant tout produit de deux $G_{\delta}$ 
denses. Un exemple de tels ensembles est le graphe d'une fonction continue et ouverte 
de domaine et d'image ouverts-ferm\'es non vides. Dans la suite, les $G_n$ et $G$ seront 
de tels graphes. Les notations et d\'efinitions qui suivent para\^\i ssent alors naturelles, 
avec le rappel qui suit. \bigskip

\noindent\bf Notation.\rm~ Soient $A$, $B$, $Z$ et $T$ des ensembles, 
$g:A\rightarrow B$ une fonction. La notation $G(g)$  d\'esignera le graphe $\mbox{Gr}(g)$ 
de $g$ si $A\times B \subseteq Z\times T$, et 
$\{(z,t)\in Z\times T~/~(t,z)\in \mbox{Gr}(g)\}$ si $A\times B \subseteq T\times Z$. On a donc 
$G(g)\subseteq Z\times T$ dans les deux cas.

\begin{defi} On dit que $(Z,T,g,(g_n))$ est une $situation\ g\acute en\acute erale$ si\smallskip

\noindent (a) $Z$ et $T$ sont des espaces polonais parfaits de dimension $0$.\smallskip

\noindent (b) $g$ et $g_n$ sont des fonctions continues et ouvertes de domaine ouvert-ferm\'e 
non vide de $Z$ et d'image ouverte-ferm\'ee de $T$, ou de domaine ouvert-ferm\'e non vide 
de $T$ et d'image ouverte-ferm\'ee de $Z$.\smallskip
 
\noindent (c) La suite $(G(g_n))$ converge vers $G(g)$.\end{defi}

 Il est d\'emontr\'e le th\'eor\`eme suivant dans [Le2] (cf th\'eor\`eme 2.3) :

\begin{thm} Soient $X$ et $Y$ des espaces polonais, $A$ un bor\'elien 
$\mbox{pot}(\borathree)$ et $\mbox{pot}(\bormthree)$ de $X\times Y$. Les conditions suivantes sont 
\'equivalentes :\smallskip

\noindent (a) Le bor\'elien $A$ n'est pas $\mbox{pot}(\bormone )$.\smallskip

\noindent (b) Il existe une situation g\'en\'erale $(Z,T,g,(g_n))$ et des injections continues 
$u : Z\rightarrow X$ et $v : T\rightarrow Y$ telles que 
$\overline{\bigcup_{n\in\omega} G(g_n)} \cap (u\times v)^{-1}(A) = 
\bigcup_{n\in\omega} G(g_n)$.\end{thm}

\vfill\eject

 Dans ce r\'esultat, la situation g\'en\'erale $(Z,T,g,(g_n))$ d\'epend de $A$, m\^eme si elle 
est toujours du m\^eme type. Dans le r\'esultat qu'on cherche \`a obtenir, annonc\'e dans 
l'introduction, le bor\'elien $A_1$ est ind\'ependant de $A$. On cherche donc 
essentiellement \`a obtenir un th\'eor\`eme d'interversion de quantificateurs, c'est-\`a-dire 
une version uniforme du th\'eor\`eme pr\'ec\'edent. Apr\`es ce rappel, il nous faut encore du 
vocabulaire.

\begin{defi} On dit que $(Z,(g_n))$ est une $bonne\ situation$ si\smallskip

\noindent (a) $Z$ est un ferm\'e parfait non vide de $\omega^\omega$.\smallskip

\noindent (b) $g_n$ est un hom\'eomorphisme de domaine et d'image ouverts-ferm\'es de $Z$. De plus, 
$\alpha <_{\mbox{lex}} g_n(\alpha )$ si $\alpha \in D_{g_n}$.\smallskip

\noindent (c) La suite $(\mbox{Gr}(g_n))$ converge vers la diagonale ${\it\Delta} (Z)$.\end{defi} 

\noindent\bf Notations.\rm ~Soit $Z\subseteq \omega^\omega$. On 
note $N_s$ l'ouvert-ferm\'e de base de $Z$ associ\'e \`a $s\in \omega^{<\omega}$ : 
$$N_s := \{\alpha\in Z~/~s\prec \alpha\}.$$
$\bullet$ Soit $(Z,(f_n))$ une bonne situation, $f_\emptyset := \mbox{Id}_{Z}$ et ${\cal R}$ la relation sur 
$\omega^{<\omega}$ definie comme suit :
$$s~{\cal R}~t~\Leftrightarrow ~\vert s\vert =\vert t\vert ~\mbox{et}~(N_s\times N_t)\cap 
\left(\mbox{Gr}(f_\emptyset)\cup \bigcup_{n\in \omega} \mbox{Gr}(f_n)\right)\not= \emptyset .$$
$\bullet$ Si $s~{\cal R}~t$, on pose 
$$m(s,t) := \mbox{min}\{m\in\omega~/~\exists~w\in \{\emptyset\}\cup \omega~~~\vert w\vert =m~~
\mbox{et}~~(N_s\times N_t)\cap \mbox{Gr}(f_w)\not= \emptyset\}\mbox{,}$$
$$n(s,t) := \mbox{min}\{n\in\omega~/~s\lceil n~{\cal R}~t\lceil n~~\mbox{et}~~m(s,t)=
m(s\lceil n,t\lceil n)\}.$$
$\bullet$ On pose $s~{\cal T}~t \Leftrightarrow s~{\cal R}~t~~\mbox{ou}~~
t~{\cal R}~s$. On dira que $c\in (\omega^{<\omega})^{<\omega}\setminus\{\emptyset\}$ 
est une ${\cal T}\mbox{-}cha\hat\imath ne$ si $\forall~i<\vert c\vert -1~~c(i)~{\cal T}~c(i+1)$.\bigskip

\noindent $\bullet$ On d\'efinit $\cal E$ comme \'etant la relation d'\'equivalence engendr\'ee par $\cal R$ :
$$s~{\cal E}~t \Leftrightarrow \exists~c~~{\cal T}\mbox{-cha\^\i ne}~~c(0)=s~~\mbox{et}~~
c(\vert c\vert -1)=t.$$

\begin{defi} On dit que $(Z,(f_n))$ est une $tr\grave es\ bonne\ situation$ si\smallskip

\noindent (a) $(Z,(f_n))$ est une bonne situation.\smallskip

\noindent (b) Si $c$ est une $\cal T$-cha\^\i ne telle que $\vert c\vert \geq 3$, $c(0)=c(\vert c\vert -1)$, et 
$c(i) \not= c(i+1)$ si $i<\vert c\vert -1$, alors il existe $i<\vert c\vert -2$ tel que $c(i) = c(i+2)$.\end{defi}

\begin{thm} Soient $(Z,(g_n))$ une bonne situation et $(\omega^\omega,(f_n))$ une 
tr\`es bonne situation. On suppose que les classes d'\'equivalence de la relation 
$\cal E$ associ\'ee \`a $(\omega^\omega,(f_n))$ sont finies. Alors il existe une fonction 
continue $u : \omega^\omega\rightarrow Z$ telle que
$$\overline{\bigcup_{n\in\omega} \mbox{Gr}(f_n)}\cap (u\times u)^{-1}
\left(\bigcup_{n\in\omega} \mbox{Gr}(g_n)\right) = \bigcup_{n\in\omega} \mbox{Gr}(f_n).$$\end{thm}

\vfill\eject

\noindent\bf D\'emonstration.\rm\ On va construire\bigskip

\noindent - Une suite $(U_s)_{s\in \omega^{<\omega}}$ d'ouverts-ferm\'es non vides de $Z$.\smallskip

\noindent - Une fonction $\Phi : \{(s,t)\in \omega^{<\omega}\times\omega^{<\omega}~/~
\vert s\vert =\vert t\vert \} \rightarrow \{\emptyset\}\cup \omega$.\bigskip

\noindent On notera, si $s~{\cal R}~t$, $w(s,t) := \Phi(s\lceil n(s,t),t\lceil n(s,t))$. On 
demande \`a ces objets de v\'erifier
$$\begin{array}{ll}
& (i)~~~U_{s^\frown i}\subseteq U_s\cr 
& (ii)~~{\delta} (U_{s^\frown i})\leq 2^{-\vert s\vert -1} \cr 
& (iii)~s~{\cal R}~t \Rightarrow \left\{\!\!\!\!\!\!
\begin{array}{ll}
 & \vert w(s,t)\vert  = m(s,t) \cr 
& U_t = g_{w(s,t)} [U_s]
\end{array}\right.
\end{array}$$
$\bullet$ Admettons ceci r\'ealis\'e. Soit $\alpha$ dans $\omega^\omega$. Comme pour $q>0$, 
${\delta} (U_{\alpha\lceil q})<2^{-q}$, $(U_{\alpha\lceil q})_q$ est une suite 
d\'ecroissante de ferm\'es non vides dont les diam\`etres tendent vers $0$. On peut 
donc d\'efinir $u : \omega^\omega\rightarrow Z$ par la formule 
$\{u(\alpha )\} = \bigcap_{q\in\omega} U_{\alpha\lceil q}$, et $u$ est continue. 
Montrons que si $(\alpha ,\beta)$ est dans $\bigcup_{n\in\omega} \mbox{Gr}(f_n)$, alors 
$(u(\alpha ),u(\beta ))$ est dans $\bigcup_{n\in\omega} \mbox{Gr}(g_n)$. Soit donc 
$n$ entier tel que  $(\alpha ,\beta)\in \mbox{Gr}(f_n)$ ; on peut trouver un entier naturel 
$m_0$ tel que ${(N_{\alpha\lceil m_0}\times N_{\beta\lceil m_0})\cap 
\mbox{Gr}(f_\emptyset) = \emptyset}$. 
Alors si $m \geq m_0$, on a $\alpha\lceil m~{\cal R}~\beta\lceil m$ et on a les \'egalit\'es 
${m(\alpha\lceil m,\beta\lceil m) = m(\alpha\lceil m_0,\beta\lceil m_0)=1}$. Posons 
$n_0 := n(\alpha\lceil m_0,\beta\lceil m_0)$. Si $p\geq n_0$, on a 
$m(\alpha\lceil p,\beta\lceil p) = m(\alpha\lceil n_0,\beta\lceil n_0)$ et 
$n(\alpha\lceil p,\beta\lceil p) = n(\alpha\lceil n_0,\beta\lceil n_0)=n_0$. Posons 
$s := \alpha\lceil n_0$ et $t := \beta\lceil n_0$. Par (iii), 
$\vert \Phi(s,t)\vert  = m(s,t) = m(\alpha\lceil m_0,\beta\lceil m_0) = 1$. On a 
$$g_{\Phi(s,t)}(u(\alpha))\in g_{\Phi(s,t)} [\bigcap_{n\geq n_0} U_{\alpha\lceil n}] 
\subseteq \bigcap_{n\geq n_0} g_{\Phi(s,t)} [U_{\alpha\lceil n}] = 
\bigcap_{n\geq n_0} U_{\beta\lceil n} = \{u(\beta)\}.$$
 D'o\`u $(u(\alpha ),u(\beta )) \in \mbox{Gr}(g_{\Phi(s,t)})$. Si $(\alpha, \beta)\in 
\overline{\bigcup_{n\in\omega} \mbox{Gr}(f_n)}\setminus (\bigcup_{n\in\omega} \mbox{Gr}(f_n))$, 
$\alpha = \beta$ et $u(\alpha )= u(\beta )$. Donc 
$(u(\alpha ), u(\beta ))\notin\bigcup_{n\in\omega} \mbox{Gr}(g_n)$.\bigskip 

\noindent $\bullet$ Montrons donc que la construction est possible. On pose 
$\Phi(\emptyset,\emptyset) := \emptyset~{\rm et}~U_\emptyset := Z$. Admettons avoir construit $U_s$ et 
$\Phi (s,t)$ pour $\vert s\vert $, $\vert t\vert \leq p$ v\'erifiant (i)-(iii), et soient $s\in \omega^p$ et 
$i\in \omega$. Posons 
$$d:\left\{\!\!
\begin{array}{ll} 
{\cal E}(s^\frown i)\times {\cal E}(s^\frown i)\!\!\!\! & 
\rightarrow \omega \cr (x,y) & \mapsto ~\mbox{min}\{\vert c\vert -1~/~c~~{\cal T}
\mbox{-cha\^\i ne,}~c(0) = x~~\mbox{et}~~c(\vert c\vert -1) = y\} 
\end{array}
\right.$$
Si $k\in\omega$, on pose $H_k := \{z\in {\cal E}(s^\frown i)~/~d(z,s^\frown i) = 
k\}$. Alors $H_k$ et le nombre de $H_k$ non vides sont finis, puisque les classes 
d'\'equivalence de $\cal E$ sont suppos\'ees finies. De 
plus, $H_k$ est non vide si $H_{k+1}$ l'est, donc on peut trouver $q$ tel que $H_0$,
...,$H_q$ soient non vides et $H_k$ soit vide si $k>q$. Posons 
$$H_k := \{z_{(k,1)},...,z_{(k,p_k)}\},~\phi : \left\{\!\!
\begin{array}{ll}
\bigcup_{k\leq q} \{ k\}\times\{1,...,p_k\}\!\!\!\! 
& \rightarrow \omega \cr 
(k,r) 
& \mapsto ({\Sigma}_{i<k}~p_i)+r  
\end{array}
\right.$$
On a donc Im$(\phi ) = \{1,...,p_0,p_0 +1,...,p_0 + p_1,...,p_0 +...+p_{q-1} +1, ..., p_0+...+p_q\}$.\bigskip
 
  On va construire par r\'ecurrence sur $n\in\{1,...,p_0+...+p_q\}$, et 
pour $k\in \{1,...,n\}$, des ouverts-ferm\'es non vides $U^n_{z_{\phi^{-1}(k)}}$ de $Z$. 
Si $z_{\phi^{-1}(k)}~{\cal R}~z_{\phi^{-1}(l)}$, on note 
$w(k,l) := w(z_{\phi^{-1}(k)},z_{\phi^{-1}(l)})$.
 
\vfill\eject

 On demande aux ouverts-ferm\'es de v\'erifier 
$$\begin{array}{ll} 
& (1)~U^n_{z_{\phi^{-1}(k)}} \subseteq U_{z_{\phi^{-1}(k)}\lceil p}\cr 
& (2)~{\delta} (U^n_{z_{\phi^{-1}(k)}})\leq 2^{-p-1}\cr
& (3)~\mbox{Si}~k,l\in\{1,...,n\}~\mbox{et}~z_{\phi^{-1}(k)}~{\cal R}~z_{\phi^{-1}(l)}\mbox{,~alors} \cr 
& ~~~~~~-\vert w(k,l)\vert  = m(z_{\phi^{-1}(k)},z_{\phi^{-1}(l)})  \cr
& ~~~~~~-U_{z_{\phi^{-1}(l)}}^n = g_{w(k,l)}[U_{z_{\phi^{-1}(k)}}^n] \cr 
& (4)~U^{n+1}_{z_{\phi^{-1}(k)}} \subseteq U^n_{z_{\phi^{-1}(k)}}~\mbox{si}~k\in\{1,...,n\}
\end{array}$$
Admettons cette construction effectu\'ee. Il restera \`a poser, si $z,~z'\in {\cal E}(s^\frown i)$, 
$U_z := U^{p_0+...+p_q}_z$, et si $\Phi(z,z')$ n'est pas encore d\'efini, on posera 
$\Phi(z,z') := \Phi(z\lceil p,z'\lceil p)$. On a 
$$U_{z_{\phi^{-1}(k)}} = U_{z_{\phi^{-1}(k)}}^{p_0+...+p_q}
\subseteq U_{z_{\phi^{-1}(k)}\lceil p}.$$
 La condition (i) est donc r\'ealis\'ee pour toute suite de ${\cal E}(s^\frown i)$. La 
condition (2) entra\^\i nera de m\^eme que (ii) est r\'ealis\'ee pour toute suite de 
${\cal E}(s^\frown i)$. Pour (iii), il suffit de remarquer que si $\tilde s~{\cal R}~\tilde t$, $\tilde s$ et $\tilde t$ sont dans la m\^eme $\cal E$-classe, et (3) donne le r\'esultat.\bigskip

\noindent $\bullet$ Montrons donc que cette nouvelle construction est possible. 
Si $n=1$, $\phi^{-1}(n)$ vaut $(0,1)$ et $z_{\phi^{-1}(n)}=s^\frown i$ ; on choisit 
pour $U^1_{z_{(0,1)}}$ un ouvert-ferm\'e non vide de $U_s$, de diam\`etre au plus 
$2^{-p-1}$. \bigskip

 Admettons avoir construit les suites finies 
$U^1_{z_{\phi^{-1}(1)}}$, ..., $U^{n-1}_{z_{\phi^{-1}(1)}}$, ..., 
$U^{n-1}_{z_{\phi^{-1}(n-1)}}$, v\'erifiant (1)-(4), ce qui est fait pour $n=2$. La suite 
$z_{\phi^{-1}(n)}$ est dans $H_{(\phi^{-1}(n))_0}$, donc on peut trouver une $\cal T$-
cha\^\i ne $c$ telle que $c(0) = s^\frown i$, $c(\vert c\vert -1) = z_{\phi^{-1}(n)}$ et 
$\vert c\vert -1 = (\phi^{-1}(n))_0$. Comme $p_0=1$, $(\phi^{-1}(n))_0\geq 1$, donc 
$\vert c\vert \geq 2$ et $c(\vert c\vert -2)\in H_{(\phi^{-1}(n))_0-1}$ ; par le choix de $\phi$, 
on peut trouver $m<n$ tel que $c(\vert c\vert -2) = z_{\phi^{-1}(m)}$. D'o\`u 
$z_{\phi^{-1}(n)}~{\cal T}~z_{\phi^{-1}(m)}$. Notons 
$$o := \left\{\!\!\!\!\!\!\!\!
\begin{array}{ll} 
& n(z_{\phi^{-1}(n)},z_{\phi^{-1}(m)})~~\mbox{si}~~z_{\phi^{-1}(n)}~{\cal R}~z_{\phi^{-1}(m)}\mbox{,}\cr 
& n(z_{\phi^{-1}(m)},z_{\phi^{-1}(n)})~~\mbox{si}~~z_{\phi^{-1}(m)}~{\cal R}~z_{\phi^{-1}(n)}.
\end{array}
\right.$$
\bf Cas 1.\rm\ $o<p+1$.\bigskip

\noindent 1.1. $z_{\phi^{-1}(m)}~{\cal R}~z_{\phi^{-1}(n)}$.\bigskip

 La suite $w(m,n) = \Phi(z_{\phi^{-1}(m)}\lceil o,z_{\phi^{-1}(n)}\lceil o)$ a d\'ej\`a 
\'et\'e d\'efinie et on a 
$$U_{z_{\phi^{-1}(n)\lceil p}} = g_{w(m,n)} [U_{z_{\phi^{-1}(m)\lceil p}}].$$
 On choisit, dans $g_{w(m,n)}[U_{z_{\phi^{-1}(m)}}^{n-1}]$, un ouvert-ferm\'e non vide 
$U_{z_{\phi^{-1}(n)}}^n$ de diam\`etre au plus $2^{-p-1}$. De 
sorte que (1), (2), et (3) pour $k = l = n$ sont r\'ealis\'ees.\bigskip

 On d\'efinit ensuite les $U_{z_{\phi^{-1}(q)}}^n$ pour $1\leq q < n$, par r\'ecurrence 
sur $d(z_{\phi^{-1}(q)},z_{\phi^{-1}(n)})$ : on choisit 
une $\cal T$-cha\^\i ne $e$ de longueur minimale telle que $e(0) = z_{\phi^{-1}(q)}$ et 
$e(\vert e\vert -1) = z_{\phi^{-1}(n)}$. Comme $\vert e\vert \geq 2$, $U_{e(1)}^n$ a \'et\'e d\'efini et 
il y a 2 cas. Soit $r$ entier compris entre $1$ et $n$ tel que 
$e(1) = z_{\phi^{-1}(r)}$. Un tel $r$ existe car la condition (b) de la d\'efinition d'une tr\`es 
bonne situation entra\^\i ne l'unicit\'e d'une $\cal T$-cha\^\i ne sans termes cons\'ecutifs 
identiques allant d'une suite \`a une autre ; cette $\cal T$-cha\^\i ne est donc de longueur 
minimale, et la d\'efinition de $\phi$ montre l'existence de $r$.

\vfill\eject

\noindent 1.1.1. $z_{\phi^{-1}(r)}~{\cal R}~z_{\phi^{-1}(q)}$.\bigskip

 On pose 
$$U_{z_{\phi^{-1}(q)}}^n := g_{w(r,q)}[U_{z_{\phi^{-1}(r)}}^n].$$
1.1.2. $z_{\phi^{-1}(q)}~{\cal R}~z_{\phi^{-1}(r)}$.\bigskip

 On pose 
$$U_{z_{\phi^{-1}(q)}}^n := g^{-1}_{w(q,r)} (U_{z_{\phi^{-1}(r)}}^n).$$
Montrons que ces d\'efinitions sont licites. On a $e(1) = z_{\phi^{-1}(r)}$, o\`u 
$1\leq r \leq n$. Si le cas $r = n$ se produit, comme $z_{\phi^{-1}(m)}$ et 
$z_{\phi^{-1}(q)}$ sont dans ${\cal E}(s^\frown i)$, l'unicit\'e de la $\cal T$-cha\^\i ne 
sans termes cons\'ecutifs identiques allant de $s^\frown i$ \`a $z_{\phi^{-1}(n)}$ 
montre que $q = m$.\bigskip

 On en d\'eduit que si $r = n$, on est dans le cas 1.1.2 puisqu'on ne peut pas avoir 
$z_{\phi^{-1}(q)}~{\cal R}~e(1)$ et $e(1)~{\cal R}~z_{\phi^{-1}(q)}$, ces deux suites \'etant 
diff\'erentes par minimalit\'e de la longueur de $v$ (si $\tilde s~{\cal R}~\tilde t$, 
on a que $\tilde s \leq_{\mbox{lex}} \tilde t$, par d\'efinition d'une tr\`es bonne 
situation).\bigskip

 Dans le cas 1.1.1, on a $r<n$ et 
$U^{n-1}_{z_{\phi^{-1}(q)}} = g_{w(r,q)}[U^{n-1}_{z_{\phi^{-1}(r)}}]$, donc 
$U^{n}_{z_{\phi^{-1}(q)}}$ est un ouvert-ferm\'e non vide de 
$U^{n-1}_{z_{\phi^{-1}(q)}}$, puisque 
$U^{n}_{z_{\phi^{-1}(r)}}\subseteq U^{n-1}_{z_{\phi^{-1}(r)}}$. De m\^eme, 
$U^{n}_{z_{\phi^{-1}(q)}}$ est un ouvert-ferm\'e non vide de 
$U^{n-1}_{z_{\phi^{-1}(q)}}$ dans le cas 1.1.2, $r<n$. Si $r=n$, $q=m$ et la m\^eme 
conclusion vaut, par le choix de $U^{n}_{z_{\phi^{-1}(n)}}$. D'o\`u la condition (4). Les 
conditions (1) et (2) pour $k=q$ en d\'ecoulent.\bigskip

 V\'erifions (3). Soient donc $k,l\leq n$ tels que $z_{\phi^{-1}(k)}~{\cal R}~z_{\phi^{-1}(l)}$, et 
$\tilde c$ (resp. $\tilde e$) la $\cal T$-cha\^\i ne ayant servi \`a d\'efinir 
$U_{z_{\phi^{-1}(k)}}^n$ (resp. $U_{z_{\phi^{-1}(l)}}^n$). On a 
$\tilde c(\vert \tilde c\vert -1) = \tilde e(\vert \tilde e\vert -1) = z_{\phi^{-1}(n)}$. 
Si $\vert \tilde c\vert  = \vert \tilde e\vert  = 1$, $k=l=n$ et (3) a \'et\'e v\'erifi\'e. Plus 
g\'en\'eralement, si $k=l$, (3) est v\'erifi\'e. Si 
$\vert \tilde c\vert  = 1$ et  $\vert \tilde e\vert  = 2$, la liaison entre $z_{\phi^{-1}(k)}$ et 
$z_{\phi^{-1}(l)}$ a d\'ej\`a \'et\'e prise en compte, par minimalit\'e des longueurs. De m\^eme 
si $\vert \tilde c\vert  = 2$ et $\vert \tilde e\vert  = 1$. Si $\vert \tilde c\vert $ et 
$\vert \tilde e\vert $ sont au moins \'egaux \`a 2, par unicit\'e de la $\cal T$-cha\^\i ne 
sans termes cons\'ecutifs identiques allant d'une suite \`a une autre, on a que 
$\tilde c(1) = \tilde e(0)$ ou $\tilde c(0) = \tilde e(1)$. L\`a 
encore, la liaison a \'et\'e prise en compte. La condition (3) est donc r\'ealis\'ee.\bigskip

\noindent 1.2. $z_{\phi^{-1}(n)}~{\cal R}~z_{\phi^{-1}(m)}$.\bigskip

 Ce cas est analogue au pr\'ec\'edent (on a 
$U_{z_{\phi^{-1}(m)\lceil p}} = g_{w(n,m)} [U_{z_{\phi^{-1}(n)\lceil p}}]$, on choisit 
$U_{z_{\phi^{-1}(n)}}^n$ dans $g_{w(n,m)}^{-1}(U_{z_{\phi^{-1}(m)}}^{n-1})$, et seul le 
cas 1.2.1 est possible si $r=n$).\bigskip

\noindent\bf Cas 2.\rm\ $o=p+1$.\bigskip

\noindent 2.1. $z_{\phi^{-1}(m)}~{\cal R}~z_{\phi^{-1}(n)}$.\bigskip

 Soit $w\in \{\emptyset\}\cup \omega$ tel que 
$(N_{z_{\phi^{-1}(m)\lceil p}}\times 
N_{z_{\phi^{-1}(n)\lceil p}})\cap \mbox{Gr}(f_w)\not=\emptyset$. On peut supposer que 
$$\vert w\vert  = m(z_{\phi^{-1}(m)\lceil p},z_{\phi^{-1}(n)\lceil p})=0\mbox{,}$$
et $w=\emptyset$, $z_{\phi^{-1}(m)\lceil p} = z_{\phi^{-1}(n)\lceil p}$. Comme  
${\mbox{Gr}(g_\emptyset) = \overline{\bigcup_{n\in\omega} \mbox{Gr}(g_{n})}\setminus 
(\bigcup_{n\in\omega} \mbox{Gr}(g_{n}))}$, on peut trouver 
$t\in \omega$ minimal tel que 
$(U_{z_{\phi^{-1}(m)}}^{n-1}\times U_{z_{\phi^{-1}(m)}}^{n-1})
\cap \mbox{Gr}(g_{t})\not=\emptyset$, et on a 
$\vert t\vert  = m(z_{\phi^{-1}(m)},z_{\phi^{-1}(n)})=1$. On pose alors 
$\Phi(z_{\phi^{-1}(m)},z_{\phi^{-1}(n)}) := t$.

\vfill\eject

 On a alors que $g_{w(m,n)} [U_{z_{\phi^{-1}(m)}}^{n-1}\cap 
g_{w(m,n)}^{-1}(U_{z_{\phi^{-1}(m)}}^{n-1})]$ est un 
ouvert-ferm\'e non vide de $U_{z_{\phi^{-1}(n)\lceil p}}$. On choisit 
$U_{z_{\phi^{-1}(n)}}^n$ dans cet ouvert-ferm\'e et on raisonne comme en 1.1.\bigskip

\noindent 2.2. $z_{\phi^{-1}(n)}~{\cal R}~z_{\phi^{-1}(m)}$.\bigskip

 On raisonne comme en 2.1, en choisissant $U_{z_{\phi^{-1}(n)}}^n$ dans 
$g_t^{-1}(U_{z_{\phi^{-1}(m)}}^{n-1})\cap U_{z_{\phi^{-1}(m)}}^{n-1}$ et en posant 
$\phi (z_{\phi^{-1}(n)},z_{\phi^{-1}(m)}) := t$.$\hfill\square$\bigskip

\noindent\bf (B) L'existence de tests.\bigskip\rm

 Nous donnons maintenant un exemple explicite, comme annonc\'e dans 
l'introduction. Nous commenons par un exemple dans $\omega^\omega\times\omega^\omega$, que nous raffinons ensuite dans un produit $Z_0\times Z_0$, o\`u $Z_0$ est plus compliqu\'e \`a 
d\'ecrire que $\omega^\omega$, mais est hom\'eomorphe \`a $2^\omega$.\bigskip

\noindent\bf Notations.\rm~Soit $(q_n)$ la suite des nombres premiers : 
$q_0 = 2$, $q_1 = 3$, $q_2 = 5$, ... On pose 
$$N : \left\{\!\!
\begin{array}{ll}
\omega^{<\omega}\!\!\!\! & \rightarrow 
\omega\cr\cr 
s 
& \mapsto \left\{\!\!\!\!\!\!
\begin{array}{ll} 
& q_0^{s(0)+1}...q_{\vert s\vert -1}^{s(\vert s\vert -1)+1}~\mbox{si}~s\not=\emptyset\mbox{,}\cr\cr 
& 0~\mbox{sinon.}
\end{array}\right.
\end{array}
\right.$$
$\bullet$ La fonction $f_\emptyset$ 
est l'identit\'e. On~ pose~ ensuite
$$f_{n}:\left\{\!\!
\begin{array}{ll}
\{\alpha\in\omega^\omega~/~\alpha (n)=1\}\!\!\!\! 
& \rightarrow\{\alpha\in\omega^\omega~/~\alpha (n) = N(\alpha\lceil n^\frown 1)\}\cr\cr  
\alpha 
& \mapsto\left\{\!\!
\begin{array}{ll}
\omega\!\!\!\! 
& \!\rightarrow\! \omega\cr 
p 
& \!\mapsto\! \left\{\!\!\!\!\!\!
\begin{array}{ll} 
& \! \alpha (p)~\mbox{si}~p\!\not=\! n\mbox{,}\cr\cr 
&\! N(\alpha \lceil (n+1))~\mbox{sinon.}
\end{array}
\right.
\end{array}
\right.
\end{array}
\right.$$
\noindent $\bullet$ On pose ensuite 
$${\cal A}_0 := \{ 1\}\mbox{,}$$ 
$${\cal A}_{n+1} := \{ 1\} \cup \left\{ N(s^\frown 1)~/~s\in \Pi_{i\leq n}~{\cal A}_i\right\}~~(n\in\omega)\mbox{,}$$
$$Z_0 := \Pi_{n\in\omega}~{\cal A}_n.$$
Alors on voit facilement par r\'ecurrence que ${\cal A}_n$ est fini et a au moins deux 
\'el\'ements si $n\geq 1$, de sorte que $Z_0$, muni de la topologie induite par celle de 
$\omega^\omega$, est hom\'eomorphe \`a $2^\omega$, comme compact m\'etrisable parfait de 
dimension 0 non vide. Il est clair que si $\alpha\in Z_0$ et $\alpha (n)=1$, alors 
$f_n (\alpha )\in Z_0$, de sorte qu'on peut remplacer $\omega^\omega$ par $Z_0$ dans la 
d\'efinition de $f_n$. On note encore $f_n$ cette nouvelle fonction, le contexte 
pr\'ecisant si on travaille dans $\omega^\omega$ ou dans $Z_0$.

\begin{thm} (1) Le couple $(\omega^\omega,(f_n))$ une tr\`es bonne situation. De 
plus, les classes d'\'equivalence de $\cal E$ sont finies.\smallskip

\noindent (2) Le couple $(Z_0,(f_n)_{n>0})$ une tr\`es bonne situation. De plus, pour tout entier 
$p$, $\Pi_{n<p}~{\cal A}_n$ est une classe pour $\cal E$.\end{thm}

\vfill\eject

\noindent\bf D\'emonstration.\rm\ Les espaces $\omega^\omega$ 
et $Z_0$ sont ferm\'es parfaits non vides de $\omega^\omega$.\bigskip

\noindent $\bullet$ Que ce soit dans $\omega^\omega$ ou $Z_0$, $f_n$ est clairement un 
hom\'eomorphisme de domaine et d'image ouverts-ferm\'es, et on a $\alpha <_{\mbox{lex}} f_n(\alpha )$ pour tout $\alpha$ de $D_{f_n}$.\bigskip

\noindent $\bullet$ Si $\alpha\!\in\! \omega^\omega$, la suite de terme g\'en\'eral 
${(\!\alpha\lceil n^\frown 1^\frown (\!\alpha (n+1),...\! ), 
\alpha\lceil n^\frown N(\alpha\lceil n^\frown 1)^\frown (\!\alpha (n+1),...\! )\! )}$ converge 
vers $(\alpha,\alpha)$, de sorte que $(\mbox{Gr}(f_n))$ converge vers ${\it\Delta} (\omega^\omega)$. De 
m\^eme si $\alpha\in Z_0$.\bigskip

\noindent $\bullet$ Montrons maintenant que pour tout entier 
$p$, $\Pi_{n<p}~{\cal A}_n$ est une classe pour $\cal E$. Il suffit de voir que 
${\cal E} (1^p) = \Pi_{n<p}~{\cal A}_n$. La condition est clairement v\'erifi\'ee pour 
$p=0$ : on a ${\cal E} (\emptyset) = \{ \emptyset\}$. Pour $p=1$, on a 
${\cal E} (1^p) = \{ 1\} = {\cal A}_0$ car on consid\`ere la suite $(f_n)_{n>0}$, de 
sorte que la premi\`ere coordonn\'ee vaut toujours 1. Soit donc $s\in {\cal E} (1^{p+1})$. 
On a bien s\^ur $s\in \Pi_{n<p+1}~{\cal A}_n$. R\'eciproquement, si $s\in \Pi_{n<p+1}~ 
{\cal A}_n$, $s\lceil p\in \Pi_{n<p}~{\cal A}_n$, donc par hypoth\`ese de r\'ecurrence, 
on peut trouver une ${\cal T}$-cha\^\i ne $v$ telle que $v(0) = 1^p$ et 
$v(\vert v\vert -1)=s\lceil p$. On a donc que $(v(i)^\frown s(p))_{i<\vert v\vert }$ est une 
${\cal T}$-cha\^\i ne, et donc que $s \in {\cal E} ({1^p}^\frown s(p))$. D'o\`u le r\'esultat 
si $s(p)=1$. Sinon, on peut trouver $t$ dans $\Pi_{n<p}~{\cal A}_n$ telle que 
$s(p) = N(t^\frown 1)$. Par hypoth\`ese de r\'ecurrence, on peut trouver une ${\cal T}$-
cha\^\i ne $w$ telle que $w(0) = 1^p$ et $w(\vert w\vert -1)=t$. Comme avant, 
$(w(i)^\frown 1)_{i<\vert w\vert }$ est une ${\cal T}$-cha\^\i ne, donc $t^\frown 1\in {\cal E} 
(1^{p+1})$. Donc $t^\frown N(t^\frown 1)\in {\cal E} (1^{p+1})$, c'est-\`a-dire 
$t^\frown s(p)\in {\cal E} (1^{p+1})$. Comme $t\in {\cal E} (1^p)$, $t^\frown s(p)\in 
{\cal E} ({1^p}^\frown s(p))$ et $s\in {\cal E} (1^{p+1})$.\bigskip

\noindent $\bullet$ Montrons maintenant les classes d'\'equivalence de $\cal E$ sont finies. Soit 
$C$ une $\cal E$-classe, $t_0\in C$ et $s_0\in C$ lexicographiquement minimale. Une 
telle suite existe car on d\'efinit, si $q<\vert t_0\vert $, 
$s_0(q)$ comme \'etant $\mbox{min}\{s(q)~/~s\in C~{\rm et}~s\lceil q = s_0\lceil q\}$. On 
montre par r\'ecurrence sur $i<\vert s_0\vert $ que 
$$\begin{array}{ll}
& (i)~~ {\cal E}(s_0\lceil (i+1))~\mbox{est~finie.}\cr 
& (ii)~\forall~s\in {\cal E}(s_0\lceil (i+1))~~s(i)\in \{s_0(i)\}\cup\{N(u^\frown s_0(i))~/~u\in {\cal E}(s_0\lceil i)\}.
\end{array}$$
Si $i=0$, ${\cal E}(s_0 (0)) = \{s_0 (0)\}$ si $s_0(0)\not=1$, 
et $\{s_0 (0)\}\cup \{N(1)\}$ sinon, d'o\`u le r\'esultat. Admettons 
ce r\'esultat pour $i<j<\vert s_0\vert $, ce qui est v\'erifi\'e pour $j=1$. 
Montrons-le pour $j$, ce qui prouvera que $C$ est finie.\bigskip

 Soit $\tilde t\in {\cal E}(s_0\lceil (j+1))$ ; il existe une 
$\cal T$-cha\^\i ne $u$ telle que $u(0)=\tilde t$ et $u(\vert u\vert -1)= s_0\lceil (j+1)$. Si 
$i<\vert u\vert -1$, comme $u(i) ~{\cal T}~ u(i+1)$, 
$u(i)\lceil j ~{\cal T}~ u(i+1)\lceil j$, donc $u(i)\lceil j ~{\cal E}~ 
u(i+1)\lceil j$ et $\tilde t\lceil j~{\cal E}~s_0\lceil j$.\bigskip

 Montrons (ii) ; (i) s'en d\'eduira car $s\lceil j \in {\cal E}(s_0\lceil j)$ 
qui est fini par hypoth\`ese de r\'ecurrence, et car $s(j)$ est dans un ensemble 
fini. On montre que
$$\forall~s,t\in {\cal E}(s_0\lceil (j+1))~~t(j)\not=s_0(j)~\mbox{ou}~
s(j)\in\{ s_0(j)\}\cup \{N(u^\frown s_0(j))~/~u\in {\cal E}(s_0\lceil j)\}.$$
On proc\`ede par r\'ecurrence sur $d(s,t)$. C'est clair pour $d(s,t) = 0$. 
Soient $s,t\in  {\cal E}(s_0\lceil (j+1))$ telles que 
$d(s,t) = k+1$. Soit $e$ une $\cal T$-cha\^\i ne telle que $e(0) = t$, 
$e(\vert e\vert -1) = s$ et $\vert e\vert  = k+2$. Soit $i<\vert e\vert $ maximal tel que 
$e(i)(j) = s_0(j)$. Si $i>0$, par hypoth\`ese de r\'ecurrence, on a le r\'esultat. On 
peut donc supposer que si $1\leq i <\vert e\vert $, $e(i)(j) \not= s_0(j)$.\bigskip

Par cons\'equent, $e(1)(j) \not= s_0(j)$ et $e(0)(j)=t(j)=s_0(j)$, donc $s_0(j)$ a 
\'et\'e modifi\'e  en $N(t\lceil j^\frown s_0(j))$, par minimalit\'e de $s_0$.

\vfill\eject

 En effet, on remarque que si $x~{\cal R}~y$, alors 
${\forall~l\in\omega\setminus\{ 0\}}$, $x^\frown l~{\cal R}~y^\frown l$. Comme $t\lceil j~{\cal E}~s_0\lceil j$, 
$$e(1)=t\lceil j^\frown e(1)(j)~{\cal E}~s_0\lceil j^\frown e(1)(j)$$ et si on pose 
$s'_0 := <s_0(j+1),...,s_0(\vert s_0\vert -1)>$, on a $e(1)^\frown s'_0~{\cal E}~s_0$ et 
$e(1)^\frown s'_0~{\cal E}~s_0\lceil j^\frown e(1)(j)^\frown s'_0$. Donc 
$e(1)(j)\geq s_0(j)$ et $e(1)(j)>s_0(j)$. Pour transformer \`a nouveau $e(1)(j)$, 
on ne peut que revenir \`a $s_0(j)$, ce qui est exclus. Donc $e(1)(j)$ reste fixe dans 
la suite et vaut $s(j)$. L'entier $s(j)$ a donc la forme voulue. D'o\`u (ii), avec 
$t = s_0\lceil (j+1)$.\bigskip

\noindent $\bullet$ Montrons maintenant que $(\omega^\omega,(f_n))$ une tr\`es bonne situation. Nous 
voulons montrer que si $c$ est une $\cal T$-cha\^\i ne telle que 
$\vert c\vert \geq 3$, $c(0)=c(\vert c\vert -1)$, et $c(i) \not= c(i+1)$ si $i<\vert c\vert -1$, alors 
il existe $i<\vert c\vert -2$ tel que $c(i) = c(i+2)$.\bigskip

 Soit $c$ un contre-exemple de longueur 
minimale, et tel que $l:= \vert c(0)\vert $ soit minimale elle aussi. Alors 
n\'ecessairement la suite 
$(c(i)(l-1))_{i<\vert c\vert }$ est non constante, et on trouve $i_1$ minimal tel que 
$c(i_1)(l-1)\not= c(i_1+1)(l-1)$ ; il y a alors deux cas.\bigskip

 Ou bien $c(i_1)(l-1)<c(i_1+1)(l-1)$, auquel cas comme on a les \'egalit\'es 
$${c(i_1)(l-1)\! =\! c(0)(l-1)\! =\! c(\vert c\vert -1)(l-1)}\mbox{,}$$
on trouve $i_2>i_1+1$ 
minimal tel que l'on ait ${c(i_1+1)(l-1)\not= c(i_2)(l-1)}$. Comme avant, on voit 
que ${c(i_1)(l-1)=c(i_2)(l-1)}$, et en fait $c(i_1)=c(i_2)$. Donc $i_1=0$ et 
$i_2=\vert c\vert -1$, par minimalit\'e de $\vert c\vert $. Par minimalit\'e encore, $\vert c\vert =3$, 
ce qui constitue la contradiction cherch\'ee (on a $c(i_1+1) = c(i_2-1)$ car il 
existe un unique entier $n$ tel que $c(i_1)^\frown 1^\omega\in 
A_{n}$, avec $n = l-1$ ; par suite, 
$c(i_1+1)^\frown 1^\omega = f_{n} (c(i_1)^\frown 1^\omega) = 
f_{n} (c(i_2)^\frown 1^\omega) = c(i_2-1)^\frown 1^\omega$).\bigskip

 Ou bien $c(i_1)(l-1)>c(i_1+1)(l-1)$, auquel cas on trouve $i_2>i_1+1$ minimal 
tel que $c(i_2)(l-1) = ... = c(\vert c\vert -1)(l-1)$. On a $c(i_1+1) = c(i_2-1)$, donc 
$c(i_1)=c(i_2)$ comme avant. D'o\`u $i_1=0$ et $i_2=\vert c\vert -1$, par minimalit\'e de 
$\vert c\vert $. Par minimalit\'e encore, $\vert c\vert =3$, ce qui constitue la contradiction 
cherch\'ee.\bigskip

\noindent $\bullet$ Il reste \`a voir que $(Z_0,(f_n)_{n>0})$ une tr\`es bonne situation pour 
achever la preuve du th\'eor\`eme. Mais ceci se voit comme pr\'ec\'edemment.$\hfill\square$
   
\begin{thm} Soit $(Z,T,g_\emptyset,(g_n))$ une situation g\'en\'erale. Alors 
il existe $u : \omega^\omega\rightarrow Z$ et $v : \omega^\omega\rightarrow T$ continues telles que 
$$\overline{\bigcup_{n\in\omega} \mbox{Gr}(f_n)}\cap (u\times v)^{-1}
\left(\bigcup_{n\in\omega} G(g_n)\right) = \bigcup_{n\in\omega} \mbox{Gr}(f_n).$$\end{thm} 

\noindent\bf D\'emonstration.\rm\ On utilisera des notations analogues \`a celles de la preuve du th\'eor\`eme 2.6, et le m\^eme sch\'ema de d\'emonstration. Les nuances sont les suivantes.

\vfill\eject

\noindent $\bullet$ On va noter $\tilde A_\emptyset$ (resp. 
$\tilde B_\emptyset$, $\tilde A_n$) le domaine de $g_\emptyset$ 
(resp. l'image de $g_\emptyset$, le domaine de $g_n$).\bigskip

\noindent $\bullet$ On va construire\bigskip

\noindent - Une suite $(U_s)_{s\in \omega^{<\omega}}$ d'ouverts non vides de $Z$, inclus dans 
$\tilde A_\emptyset$ ou $\tilde B_\emptyset$.\smallskip

\noindent - Une suite $(V_s)_{s\in \omega^{<\omega}}$ d'ouverts non vides de $T$, inclus dans 
$\tilde A_\emptyset$ ou $\tilde B_\emptyset$.\smallskip

\noindent - Une fonction $\Phi : \{(s,t)\in \omega^{<\omega}\times\omega^{<\omega}~/~
\vert s\vert =\vert t\vert \} \rightarrow \{\emptyset\}\cup \omega$.\bigskip

 On demande \`a ces objets de v\'erifier
$$\begin{array}{ll}
& (i)~~~\overline{U_{s^\frown i}}\times \overline{V_{s^\frown i}}
\subseteq U_s\times V_s\cr 
& (ii)~~{\delta} (U_{s^\frown i}), {\delta} (V_{s^\frown i}) \leq 2^{-\vert s\vert -1} \cr 
& (iii)~s~{\cal R}~t \Rightarrow \left\{\!\!\!\!\!\!
\begin{array}{ll}
 & \vert w(s,t)\vert  = m(s,t) \cr 
 & V_t = g_{w(s,t)} [U_s]~\mbox{si}~\tilde A_{w(s,t)}\subseteq Z\cr 
 & U_s = g_{w(s,t)} [V_t]~\mbox{si}~\tilde A_{w(s,t)}\subseteq T 
\end{array}
\right.
\end{array}$$
$\bullet$ Admettons ceci r\'ealis\'e. On d\'efinit $u : \omega^\omega\rightarrow Z$ et 
$v : \omega^\omega\rightarrow T$ par les formules 
$$\{u(\alpha )\} = \bigcap_{q\in\omega} U_{\alpha\lceil q}\mbox{,}$$
$\{v(\alpha )\} = \bigcap_{q\in\omega} V_{\alpha\lceil q}$. Montrons que si 
$(\alpha ,\beta)\in\bigcup_{n\in\omega} \mbox{Gr}(f_n)$ (resp. $\mbox{Gr}(f_\emptyset)$), alors 
$(u(\alpha ),v(\beta ))$ est dans $\bigcup_{n\in\omega} G(g_n)$ (resp. $G(g_\emptyset)$). Soit donc 
$w$ dans $\{\emptyset\}\cup\omega$ tel que  $(\alpha ,\beta)\in \mbox{Gr}(f_w)$ ; on peut 
trouver un entier naturel $m_0$ tel que ${(N_{\alpha\lceil m_0}\times N_{\beta\lceil 
m_0})\cap \bigcup_{s\in\{\emptyset\}\cup\omega, \vert s\vert <\vert w\vert } \mbox{Gr}(f_s) = \emptyset}$. 
Alors si $m \geq m_0$, on a $\alpha\lceil m~{\cal R}~\beta\lceil m$ et on a l'\'egalit\'e 
${m(\alpha\lceil m,\beta\lceil m) = m(\alpha\lceil m_0,\beta\lceil m_0)}=\vert w\vert $. Par (iii), on a 
$$\vert \Phi(s,t)\vert  = m(s,t) = m(\alpha\lceil m_0,\beta\lceil m_0) = \vert w\vert  = 
\left\{\!\!\!\!\!\!\!\!
\begin{array}{ll}
& 0~\mbox{si}~\beta = f_\emptyset (\alpha )\mbox{,}\cr 
& 1~\mbox{si}~\exists~n~~\beta = f_n (\alpha ).
\end{array}
\right.$$ 
Si $\tilde A_{\Phi(s,t)}\subseteq Z$, on a 
$$g_{\Phi(s,t)}(u(\alpha))\in g_{\Phi(s,t)} [\bigcap_{n\geq n_0} U_{\alpha\lceil n}] 
\subseteq \bigcap_{n\geq n_0} g_{\Phi(s,t)} [U_{\alpha\lceil n}] = 
\bigcap_{n\geq n_0} V_{\beta\lceil n} = \{v(\beta)\}.$$
Si $\tilde A_{\Phi(s,t)}\subseteq T$, on a 
$$g_{\Phi(s,t)}(v(\beta))\in g_{\Phi(s,t)} [\bigcap_{n\geq n_0} V_{\beta\lceil n}] 
\subseteq \bigcap_{n\geq n_0} g_{\Phi(s,t)} [V_{\beta\lceil n}] = 
\bigcap_{n\geq n_0} U_{\alpha\lceil n} = \{u(\alpha)\}.$$
D'o\`u $(u(\alpha ),v(\beta )) \in G(g_{\Phi(s,t)})$.\bigskip
 
\noindent $\bullet$ Montrons donc que la construction est possible. On pose 
$$\Phi(\emptyset,\emptyset) := \emptyset~~\mbox{et}~~(U_\emptyset, V_\emptyset) := 
\left\{\!\!\!\!\!\!\!\!
\begin{array}{ll} 
& (\tilde A_\emptyset, \tilde B_\emptyset)~\mbox{si}~\tilde A_\emptyset\subseteq Z\mbox{,}\cr 
& (\tilde B_\emptyset, \tilde A_\emptyset)~\mbox{si}~\tilde A_\emptyset\subseteq T.
\end{array}
\right.$$

\vfill\eject

 Par le th\'eor\`eme 2.7, les classes d'\'equivalence de $\cal E$ sont finies. On peut donc 
d\'efinir $H_k$ et $\phi$ comme dans la preuve du th\'eor\`eme 2.6. On va construire par r\'ecurrence sur $n\in\{1,...,p_0+...+p_q\}$, et 
pour $k\in \{1,...,n\}$, des ouverts non vides $U^n_{z_{\phi^{-1}(k)}}$ 
(resp. $V^n_{z_{\phi^{-1}(k)}}$) de $Z$ (resp. $T$). On demande aux ouverts de v\'erifier 
$$\begin{array}{ll} 
& (1)~\overline{U^n_{z_{\phi^{-1}(k)}}}\times \overline{V^n_{z_{\phi^{-1}(k)}}}
\subseteq U_{z_{\phi^{-1}(k)}\lceil p}\times V_{z_{\phi^{-1}(k)}\lceil p}\cr 
& (2)~{\delta} (U^n_{z_{\phi^{-1}(k)}}), {\delta} (V^n_{z_{\phi^{-1}(k)}})\leq 2^{-p-1}\cr
& (3)~\mbox{Si}~k,l\in\{1,...,n\}~\mbox{et}~z_{\phi^{-1}(k)}~{\cal R}~z_{\phi^{-1}(l)}
\mbox{,~alors} \cr 
& ~~~~~~-\vert w(k,l)\vert  = m(z_{\phi^{-1}(k)},z_{\phi^{-1}(l)})  \cr
& ~~~~~~-\mbox{si}~\tilde A_{w(k,l)}\subseteq Z\mbox{,~alors}~
V_{z_{\phi^{-1}(l)}}^n = g_{w(k,l)}[U_{z_{\phi^{-1}(k)}}^n] \cr 
& ~~~~~~-\mbox{si}~\tilde A_{w(k,l)}\subseteq T\mbox{,~alors}~
U_{z_{\phi^{-1}(k)}}^n = g_{w(k,l)}[V_{z_{\phi^{-1}(l)}}^n] \cr 
& (4)~U^{n+1}_{z_{\phi^{-1}(k)}} \times V^{n+1}_{z_{\phi^{-1}(k)}} \subseteq 
U^n_{z_{\phi^{-1}(k)}}\times V^n_{z_{\phi^{-1}(k)}}~\mbox{si}~k\in\{1,...,n\}
\end{array}$$
$\bullet$ Montrons donc que cette nouvelle construction est possible. 
Si $n=1$, on choisit pour $U^1_{z_{(0,1)}}$ un ouvert non vide de $\tilde A_\emptyset$, 
de diam\`etre au plus $2^{-p-1}$, tel que 
${\delta} (g_\emptyset[U^1_{z_{(0,1)}}])\leq 2^{-p-1}$, $\overline{
g_\emptyset [U^1_{z_{(0,1)}}]} \subseteq V_s$ et 
$\overline{U^1_{z_{(0,1)}}}\subseteq U_s$, et on pose $V^1_{z_{(0,1)}} := 
g_\emptyset [U^1_{z_{(0,1)}}]$. Ceci si $\tilde A_\emptyset \subseteq Z$. Si 
$\tilde A_\emptyset \subseteq T$, on choisit pour $V^1_{z_{(0,1)}}$ un ouvert non vide 
de $\tilde A_\emptyset$, de diam\`etre au plus $2^{-p-1}$, tel que ${\delta} (
g_\emptyset [V^1_{z_{(0,1)}}])\leq 2^{-p-1}$, $\overline{
g_\emptyset [V^1_{z_{(0,1)}}]} \subseteq U_s$ et 
$\overline{V^1_{z_{(0,1)}}}\subseteq V_s$, et on pose $U^1_{z_{(0,1)}} := 
g_\emptyset [V^1_{z_{(0,1)}}]$.\bigskip

 Admettons avoir construit les suites finies 
$U^1_{z_{\phi^{-1}(1)}}$, $V^1_{z_{\phi^{-1}(1)}}$, ..., $U^{n-1}_{z_{\phi^{-1}(1)}}$, 
$V^{n-1}_{z_{\phi^{-1}(1)}}$, ..., $U^{n-1}_{z_{\phi^{-1}(n-1)}}$, 
$V^{n-1}_{z_{\phi^{-1}(n-1)}}$ v\'erifiant (1)-(4).\bigskip

\noindent\bf Cas 1.\rm\ $o<p+1$.\bigskip

\noindent 1.1. $z_{\phi^{-1}(m)}~{\cal R}~z_{\phi^{-1}(n)}$ et 
$\tilde A_{w(m,n)}\subseteq Z$.\bigskip

 La suite $w(m,n) = \Phi(z_{\phi^{-1}(m)}\lceil o,z_{\phi^{-1}(n)}\lceil o)$ a d\'ej\`a 
\'et\'e d\'efinie et on a 
$$V_{z_{\phi^{-1}(n)\lceil p}} = g_{w(m,n)} [U_{z_{\phi^{-1}(m)\lceil p}}].$$
\noindent 1.1.1. $\tilde A_\emptyset\subseteq Z$.\bigskip

 On choisit, dans $U_{z_{\phi^{-1}(n)\lceil p}} \cap g_\emptyset^{-1}(g_{w(m,n)}
[U_{z_{\phi^{-1}(m)}}^{n-1}])$, un ouvert non vide $U_{z_{\phi^{-1}(n)}}^n$ tel que 
$$\overline{U_{z_{\phi^{-1}(n)}}^n}\times \overline{g_\emptyset
[U_{z_{\phi^{-1}(n)}}^n]}\subseteq U_{z_{\phi^{-1}(n)\lceil p}}\times 
V_{z_{\phi^{-1}(n)\lceil p}}\mbox{,}$$
${\delta} (U_{z_{\phi^{-1}(n)}}^n)\leq 2^{-p-1}$ et \'egalement 
${{\delta} (g_\emptyset[U_{z_{\phi^{-1}(n)}}^n])\leq 2^{-p-1}}$. On pose 
$V_{z_{\phi^{-1}(n)}}^n := g_\emptyset[U_{z_{\phi^{-1}(n)}}^n]$, de sorte que 
(1), (2), et (3) pour $k = l = n$ sont r\'ealis\'ees.\bigskip

 On d\'efinit ensuite les $U_{z_{\phi^{-1}(q)}}^n$ et $V_{z_{\phi^{-1}(q)}}^n$ pour 
$1\leq q < n$, par r\'ecurrence sur $d(z_{\phi^{-1}(q)},z_{\phi^{-1}(n)})$. On choisit 
une $\cal T$-cha\^\i ne $e$ de longueur minimale telle que $e(0) = z_{\phi^{-1}(q)}$ et 
$e(\vert e\vert -1) = z_{\phi^{-1}(n)}$. Comme $\vert e\vert \geq 2$, $U_{e(1)}^n$ et $V_{e(1)}^n$ 
ont \'et\'e d\'efinis et il y a 4 cas. Soit $r$ entier compris entre $1$ et $n$ tel que 
$e(1) = z_{\phi^{-1}(r)}$.

\vfill\eject

\noindent 1.1.1.1. $z_{\phi^{-1}(r)}~{\cal R}~z_{\phi^{-1}(q)}$ et $\tilde A_{w(r,q)}\subseteq Z$.\bigskip

 On pose 
$$\left\{\!\!\!\!\!\!
\begin{array}{ll} 
& V_{z_{\phi^{-1}(q)}}^n := g_{w(r,q)}[U_{z_{\phi^{-1}(r)}}^n]\mbox{,}\cr 
& U_{z_{\phi^{-1}(q)}}^n := U_{z_{\phi^{-1}(q)}}^{n-1}\cap g^{-1}_\emptyset 
(V_{z_{\phi^{-1}(q)}}^n).
\end{array}
\right.$$
1.1.1.2. $z_{\phi^{-1}(r)}~{\cal R}~z_{\phi^{-1}(q)}$ et $\tilde A_{w(r,q)}\subseteq T$.\bigskip

 On pose 
$$\left\{\!\!\!\!\!\!
\begin{array}{ll} 
& V_{z_{\phi^{-1}(q)}}^n := V_{z_{\phi^{-1}(q)}}^{n-1}\cap g^{-1}_{w(r,q)} 
(U_{z_{\phi^{-1}(r)}}^n)\mbox{,}\cr 
& U_{z_{\phi^{-1}(q)}}^n := U_{z_{\phi^{-1}(q)}}^{n-1}\cap g^{-1}_\emptyset 
(V_{z_{\phi^{-1}(q)}}^n).
\end{array}
\right.$$
1.1.1.3. $z_{\phi^{-1}(q)}~{\cal R}~z_{\phi^{-1}(r)}$ et $\tilde A_{w(q,r)}\subseteq Z$.\bigskip

 On pose 
$$\left\{\!\!\!\!\!\!
\begin{array}{ll} 
& U_{z_{\phi^{-1}(q)}}^n := U_{z_{\phi^{-1}(q)}}^{n-1}\cap g^{-1}_{w(q,r)} 
(V_{z_{\phi^{-1}(r)}}^n)\mbox{,}\cr 
& V_{z_{\phi^{-1}(q)}}^n := g_\emptyset[U_{z_{\phi^{-1}(q)}}^n].
\end{array}
\right.$$
1.1.1.4. $z_{\phi^{-1}(q)}~{\cal R}~z_{\phi^{-1}(r)}$ et $\tilde A_{w(q,r)}\subseteq T$.\bigskip

 On pose 
$$\left\{\!\!\!\!\!\!
\begin{array}{ll} 
& U_{z_{\phi^{-1}(q)}}^n := g_{w(q,r)}[V_{z_{\phi^{-1}(r)}}^n]\mbox{,}\cr 
& V_{z_{\phi^{-1}(q)}}^n := g_\emptyset[U_{z_{\phi^{-1}(q)}}^n].
\end{array}
\right.$$
Montrons que ces d\'efinitions sont licites. Si $r = n$, on est dans le cas 1.1.1.3.\bigskip

 Dans le cas 1.1.1.1, $U_{e(1)}^n\subseteq U_{e(1)}^{n-1}~\mbox{et}~
V_{z_{\phi^{-1}(q)}}^{n-1} = g_{w(r,q)} [U_{e(1)}^{n-1}]$. Par suite, la 
d\'efinition de $V_{z_{\phi^{-1}(q)}}^n$ est licite, et c'est un ouvert 
non vide de $V_{z_{\phi^{-1}(q)}}^{n-1}$. On a 
$z_{\phi^{-1}(q)}~{\cal R}~z_{\phi^{-1}(q)}$, 
$n(z_{\phi^{-1}(q)},z_{\phi^{-1}(q)}) = 0$ et $\tilde A_\emptyset 
\subseteq Z$. On en d\'eduit que $V_{z_{\phi^{-1}(q)}}^{n-1} = g_\emptyset 
[U_{z_{\phi^{-1}(q)}}^{n-1}]$, et que $U_{z_{\phi^{-1}(q)}}^n$ est un 
ouvert non vide de $U_{z_{\phi^{-1}(q)}}^{n-1}$, bien d\'efini.\bigskip

 Dans le cas 1.1.1.2, la d\'efinition de $V_{z_{\phi^{-1}(q)}}^n$ est licite, et c'est 
un ouvert de $V_{z_{\phi^{-1}(q)}}^{n-1}$. De m\^eme, la d\'efinition de 
$U_{z_{\phi^{-1}(q)}}^n$ est licite, et c'est un ouvert de 
$U_{z_{\phi^{-1}(q)}}^{n-1}$. Comme dans le cas pr\'ec\'edent, la non-vacuit\'e de 
$V_{z_{\phi^{-1}(q)}}^n$ entra\^\i ne celle de $U_{z_{\phi^{-1}(q)}}^n$. Celle de 
$V_{z_{\phi^{-1}(q)}}^n$ r\'esulte du fait que 
$$U_{z_{\phi^{-1}(r)}}^{n-1} = g_{w(r,q)}[V_{z_{\phi^{-1}(q)} }^{n-1}]$$ 
et $U_{z_{\phi^{-1}(r)}}^n\subseteq U_{z_{\phi^{-1}(r)}}^{n-1}$.\bigskip

 Dans le cas 1.1.1.3, la d\'efinition de $U_{z_{\phi^{-1}(q)}}^n$ est licite, et c'est 
un ouvert de $U_{z_{\phi^{-1}(q)}}^{n-1}$, donc de $U_{z_{\phi^{-1}(q)}\lceil p}$ et 
de $\tilde A_\emptyset$. Par suite, la d\'efinition de $V_{z_{\phi^{-1}(q)}}^n$ est 
licite, et c'est un ouvert de $V_{z_{\phi^{-1}(q)}}^{n-1}$. La non-vacuit\'e de 
$U_{z_{\phi^{-1}(q)}}^n$ entra\^\i ne celle de $V_{z_{\phi^{-1}(q)}}^n$. Celle de 
$U_{z_{\phi^{-1}(q)}}^n$ r\'esulte du fait que 
$$V_{z_{\phi^{-1}(r)}}^{n-1} = g_{w(q,r)}[U_{z_{\phi^{-1}(q)} }^{n-1}]$$ 
et $V_{z_{\phi^{-1}(r)}}^n\subseteq V_{z_{\phi^{-1}(r)}}^{n-1}$, si $r<n$.

\vfill\eject

 Si $r=n$, on a $U_{z_{\phi^{-1}(r)}}^n\subseteq U_{z_{\phi^{-1}(r)}\lceil p}\cap 
g_\emptyset^{-1}(g_{w(q,r)}[U_{z_{\phi^{-1}(q)}}^{n-1}])$, par le choix de  
$U_{z_{\phi^{-1}(n)}}^n$, puisque $q=m$. D'o\`u $V_{z_{\phi^{-1}(r)}}^n\subseteq 
g_{w(q,r)}[U_{z_{\phi^{-1}(q)}}^{n-1}]$ et $U_{z_{\phi^{-1}(q)}}^n$ est non 
vide.\bigskip

Dans le cas 1.1.1.4, on a $V_{e(1)}^n\subseteq V_{e(1)}^{n-1}~\mbox{et}~
U_{z_{\phi^{-1}(q)}}^{n-1} = g_{w(q,r)} [V_{e(1)}^{n-1}]$. Par suite, la 
d\'efinition de $U_{z_{\phi^{-1}(q)}}^n$ est licite, et c'est un ouvert 
non vide de $U_{z_{\phi^{-1}(q)}}^{n-1}$. On a comme en 1.1.1.1 que 
$$V_{z_{\phi^{-1}(q)}}^{n-1} = g_\emptyset [U_{z_{\phi^{-1}(q)}}^{n-1}]\mbox{,}$$ 
et que $V_{z_{\phi^{-1}(q)}}^n$ est un ouvert non vide de $V_{z_{\phi^{-1}(q)}}^{n-1}$, 
bien d\'efini.\bigskip

 On a donc montr\'e que dans les 4 cas, $U_{z_{\phi^{-1}(q)}}^n$ (resp. 
$V_{z_{\phi^{-1}(q)}}^n$) est un ouvert non vide de $U_{z_{\phi^{-1}(q)}}^{n-1}$ 
(resp. $V_{z_{\phi^{-1}(q)}}^{n-1}$), bien d\'efini. D'o\`u (4). V\'erifions (3) ; si $k=l$, 
(3) est v\'erifi\'e : c'est 
clair dans les cas 1.1.1.3 et 1.1.1.4, et dans les deux autres cas, on a clairement 
$g_{w(k,k)} [U_{z_{\phi^{-1}(k)}}^n]\subseteq  V_{z_{\phi^{-1}(k)}}^n$ ; si 
$y\in V_{z_{\phi^{-1}(k)}}^n$, $y\in V_{z_{\phi^{-1}(k)}}^{n-1}$ et il existe $x\in 
U_{z_{\phi^{-1}(k)}}^{n-1}$ tel que $y = g_\emptyset (x)$. Par suite, 
$x\in U_{z_{\phi^{-1}(k)}}^n$ et $V_{z_{\phi^{-1}(k)}}^n = g_\emptyset 
[U_{z_{\phi^{-1}(k)}}^n] = g_{w(k,k)} [U_{z_{\phi^{-1}(k)}}^n]$. La condition (3) est 
donc r\'ealis\'ee (c'est clair 
dans les cas 1.1.1.1 et 1.1.1.4 ; dans le cas 1.1.1.2, on utilise le fait que 
$U_{z_{\phi^{-1}(r)}}^n$ est inclus dans $U_{z_{\phi^{-1}(r)}}^{n-1}$ ; dans le cas 
1.1.1.3, on utilise le fait que 
$V_{z_{\phi^{-1}(r)}}^n$ est inclus dans $V_{z_{\phi^{-1}(r)}}^{n-1}$ si $r<n$ ; si 
$r=n$, on utilise le fait que $V_{z_{\phi^{-1}(r)}}^n$ est inclus dans 
$g_{w(q,r)}[U_{z_{\phi^{-1}(q)}}^{n-1}]$).\bigskip

\noindent 1.1.2. $\tilde A_\emptyset\subseteq T$.\bigskip

On a $g_\emptyset[g_{w(m,n)}[U_{z_{\phi^{-1}(m)}}^{n-1}]]\subseteq 
U_{z_{\phi^{-1}(n)}\lceil p}$, donc $g_{w(m,n)}[U_{z_{\phi^{-1}(m)}}^{n-1}]
\cap g_\emptyset^{-1}(U_{z_{\phi^{-1}(n)}\lceil p})$ est non vide. On choisit, 
 dans cet ouvert non vide, un ouvert $V_{z_{\phi^{-1}(n)}}^n$ non vide tel que l'on 
ait l'inclusion 
${\overline{g_\emptyset[V_{z_{\phi^{-1}(n)}}^n]}\times 
\overline{V_{z_{\phi^{-1}(n)}}^n} 
\subseteq U_{z_{\phi^{-1}(n)\lceil p}}\times 
V_{z_{\phi^{-1}(n)\lceil p}},~{\delta} (V_{z_{\phi^{-1}(n)}}^n)\leq 2^{-p-1}}$ et 
$${{\delta} (g_\emptyset[V_{z_{\phi^{-1}(n)}}^n])\leq 2^{-p-1}}.$$ 
On pose 
$U_{z_{\phi^{-1}(n)}}^n := g_\emptyset[V_{z_{\phi^{-1}(n)}}^n]$, de sorte que 
(1), (2), et (3) pour $k = l = n$ sont r\'ealis\'ees.\bigskip

 On proc\`ede alors de mani\`ere analogue \`a celle du cas 1.1.1 ; l\`a encore, on trouve 4 
cas.\bigskip

\noindent 1.1.2.1. $z_{\phi^{-1}(r)}~{\cal R}~z_{\phi^{-1}(q)}$ et $\tilde A_{w(r,q)}\subseteq Z$.\bigskip

 On pose 
$$\left\{\!\!\!\!\!\!
\begin{array}{ll} 
& V_{z_{\phi^{-1}(q)}}^n := g_{w(r,q)}[U_{z_{\phi^{-1}(r)}}^n]\mbox{,}\cr 
& U_{z_{\phi^{-1}(q)}}^n := g_\emptyset[V_{z_{\phi^{-1}(q)}}^n].
\end{array}\right.$$
1.1.2.2. $z_{\phi^{-1}(r)}~{\cal R}~z_{\phi^{-1}(q)}$ et $\tilde 
A_{w(r,q)}\subseteq T$.\bigskip

 On pose 
$$\left\{\!\!\!\!\!\!
\begin{array}{ll} 
& V_{z_{\phi^{-1}(q)}}^n := V_{z_{\phi^{-1}(q)}}^{n-1}\cap g^{-1}_{w(r,q)} 
(U_{z_{\phi^{-1}(r)}}^n)\mbox{,}\cr 
& U_{z_{\phi^{-1}(q)}}^n := g_\emptyset[V_{z_{\phi^{-1}(q)}}^n].
\end{array}\right.$$
1.1.2.3. $z_{\phi^{-1}(q)}~{\cal R}~z_{\phi^{-1}(r)}$ et $\tilde 
A_{w(q,r)}\subseteq Z$.\bigskip

 On pose 
$$\left\{\!\!\!\!\!\!
\begin{array}{ll} 
& U_{z_{\phi^{-1}(q)}}^n := U_{z_{\phi^{-1}(q)}}^{n-1}\cap g^{-1}_{w(q,r)} 
(V_{z_{\phi^{-1}(r)}}^n)\mbox{,}\cr 
& V_{z_{\phi^{-1}(q)}}^n := V_{z_{\phi^{-1}(q)}}^{n-1}\cap g^{-1}_\emptyset 
(U_{z_{\phi^{-1}(q)}}^n).
\end{array}
\right.$$
1.1.2.4. $z_{\phi^{-1}(q)}~{\cal R}~z_{\phi^{-1}(r)}$ et $\tilde 
A_{w(q,r)}\subseteq T$.\bigskip

 On pose 
$$\left\{\!\!\!\!\!\!
\begin{array}{ll} 
& U_{z_{\phi^{-1}(q)}}^n := g_{w(q,r)}[V_{z_{\phi^{-1}(r)}}^n]\mbox{,}\cr 
& V_{z_{\phi^{-1}(q)}}^n := V_{z_{\phi^{-1}(q)}}^{n-1}\cap g^{-1}_\emptyset 
(U_{z_{\phi^{-1}(q)}}^n). 
\end{array}
\right.$$
La v\'erification des conditions (1) \`a (4) est alors analogue \`a celle du cas 1.1.1 ; pour 
v\'erifier (3) dans le cas 1.1.2.3, $r=n$, on utilise le fait que 
$V_{z_{\phi^{-1}(r)}}^n\subseteq g_{w(q,r)}[U_{z_{\phi^{-1}(q)}}^{n-1}]$.\bigskip

\noindent 1.2. $z_{\phi^{-1}(m)}~{\cal R}~z_{\phi^{-1}(n)}$ et 
$\tilde A_{w(m,n)}\subseteq T$.\bigskip

 La suite $w(m,n) = \Phi(z_{\phi^{-1}(m)}\lceil o,z_{\phi^{-1}(n)}\lceil o)$ a d\'ej\`a 
\'et\'e d\'efinie et on a 
$$U_{z_{\phi^{-1}(m)\lceil p}} = g_{w(m,n)} [V_{z_{\phi^{-1}(n)\lceil p}}].$$
\noindent 1.2.1. $\tilde A_\emptyset\subseteq Z$.\bigskip

 On choisit, dans $U_{z_{\phi^{-1}(n)\lceil p}} \cap g_\emptyset^{-1}(g_{w(m,n)}^{-1}
(U_{z_{\phi^{-1}(m)}}^{n-1}))$, un ouvert non vide $U_{z_{\phi^{-1}(n)}}^n$ tel que 
$$\overline{U_{z_{\phi^{-1}(n)}}^n}\times \overline{g_\emptyset
[U_{z_{\phi^{-1}(n)}}^n]}\subseteq U_{z_{\phi^{-1}(n)\lceil p}}\times 
V_{z_{\phi^{-1}(n)\lceil p}}\mbox{,}$$
${\delta} (U_{z_{\phi^{-1}(n)}}^n)\leq 2^{-p-1}$ et \'egalement 
${{\delta} (g_\emptyset[U_{z_{\phi^{-1}(n)}}^n])\leq 2^{-p-1}}$. On pose 
$V_{z_{\phi^{-1}(n)}}^n := g_\emptyset[U_{z_{\phi^{-1}(n)}}^n]$, de sorte que 
(1), (2), et (3) pour $k = l = n$ sont r\'ealis\'ees.  On proc\`ede alors de mani\`ere analogue 
\`a celle du cas 1.1.1 ; on trouve les 4 m\^emes cas qu'en 1.1.1, et on adopte les m\^emes 
d\'efinitions. Cette fois-ci, la diff\'erence est que si $r=n$, on est dans le cas 
1.2.1.4. La v\'erification des conditions (1) \`a (4) est alors analogue \`a celle du cas 
1.1.1 ; pour v\'erifier que 
$U_{z_{\phi^{-1}(q)}}^n\subseteq U_{z_{\phi^{-1}(q)}}^{n-1}$ dans le cas 1.2.1.4, 
$r=n$, on utilise le fait que 
$g_{w(q,r)}[g_\emptyset [U_{z_{\phi^{-1}(n)}}^n]]\subseteq 
U_{z_{\phi^{-1}(q)}}^{n-1}$.\bigskip

\noindent 1.2.2. $\tilde A_\emptyset\subseteq T$.\bigskip

On a que $g_{w(m,n)}^{-1}(U_{z_{\phi^{-1}(m)}}^{n-1})
\cap g_\emptyset^{-1}(U_{z_{\phi^{-1}(n)}\lceil p})$ est non vide. On proc\`ede 
alors comme dans le cas 1.1.2 ; on a, dans le cas 1.2.2.4, $r=n$, 
$U_{z_{\phi^{-1}(q)}}^n\subseteq U_{z_{\phi^{-1}(q)}}^{n-1}$ puisqu'on a 
$V_{z_{\phi^{-1}(n)}}^n\subseteq g_{w(q,r)}^{-1}(U_{z_{\phi^{-1}(q)}}^{n-1})$.\bigskip

\noindent 1.3. $z_{\phi^{-1}(n)}~{\cal R}~z_{\phi^{-1}(m)}$ et 
$\tilde A_{w(n,m)}\subseteq Z$.\bigskip

 La suite $w(n,m) = \Phi(z_{\phi^{-1}(n)}\lceil o,z_{\phi^{-1}(m)}\lceil o)$ a d\'ej\`a 
\'et\'e d\'efinie et on a 
$$V_{z_{\phi^{-1}(m)\lceil p}} = g_{w(n,m)} [U_{z_{\phi^{-1}(n)\lceil p}}].$$ 

\vfill\eject

\noindent 1.3.1. $\tilde A_\emptyset\subseteq Z$.\bigskip

 On choisit, dans $U_{z_{\phi^{-1}(n)\lceil p}} \cap g_{w(n,m)}^{-1}(
V_{z_{\phi^{-1}(m)}}^{n-1})$, un ouvert non vide $U_{z_{\phi^{-1}(n)}}^n$ comme 
avant.  On proc\`ede alors de mani\`ere analogue 
\`a celle du cas 1.1.1 ; on trouve les 4 m\^emes cas qu'en 1.1.1, et on adopte les m\^emes 
d\'efinitions. Cette fois-ci, la diff\'erence est que si $r=n$, on est dans le cas 
1.3.1.1. La v\'erification des conditions (1) \`a (4) est alors analogue \`a celle du cas 
1.1.1 ; pour v\'erifier que 
$V_{z_{\phi^{-1}(q)}}^n\subseteq V_{z_{\phi^{-1}(q)}}^{n-1}$ dans le cas 1.3.1.1, 
$r=n$, on utilise le fait que 
$U_{z_{\phi^{-1}(r)}}^n\subseteq g_{w(r,q)}^{-1}(
V_{z_{\phi^{-1}(q)}}^{n-1})$.\bigskip

\noindent 1.3.2. $\tilde A_\emptyset\subseteq T$.\bigskip

On a que $g_\emptyset^{-1}(g_{w(n,m)}^{-1}(V_{z_{\phi^{-1}(m)}}^{n-1}))
\cap V_{z_{\phi^{-1}(n)}\lceil p}$ est non vide. On proc\`ede 
alors comme dans le cas 1.1.2 ; on a, dans le cas 1.3.2.1, $r=n$, 
$V_{z_{\phi^{-1}(q)}}^n\subseteq V_{z_{\phi^{-1}(q)}}^{n-1}$ puisqu'on a   
$U_{z_{\phi^{-1}(n)}}^n\subseteq g_{w(r,q)}^{-1}(V_{z_{\phi^{-1}(q)}}^{n-1})$.\bigskip

\noindent 1.4. $z_{\phi^{-1}(n)}~{\cal R}~z_{\phi^{-1}(m)}$ et 
$\tilde A_{w(n,m)}\subseteq T$.\bigskip

 La suite $w(n,m) = \Phi(z_{\phi^{-1}(n)}\lceil o,z_{\phi^{-1}(m)}\lceil o)$ a d\'ej\`a 
\'et\'e d\'efinie et on a 
$$U_{z_{\phi^{-1}(n)\lceil p}} = g_{w(n,m)} [V_{z_{\phi^{-1}(m)\lceil p}}].$$ 
\noindent 1.4.1. $\tilde A_\emptyset\subseteq Z$.\bigskip

 On choisit, dans $g_{w(n,m)}[V_{z_{\phi^{-1}(m)}}^{n-1}]$, un ouvert non vide 
$U_{z_{\phi^{-1}(n)}}^n$ comme avant.  On proc\`ede alors de mani\`ere analogue 
\`a celle du cas 1.1.1 ; on trouve les 4 m\^emes cas qu'en 1.1.1, et on adopte les m\^emes 
d\'efinitions. Cette fois-ci, la diff\'erence est que si $r=n$, on est dans le cas 
1.4.1.2. La v\'erification des conditions (1) \`a (4) est alors analogue \`a celle du cas 
1.1.1 ; pour v\'erifier que la condition (3) est v\'erifi\'ee dans le cas 1.4.1.2, 
$r=n$, on utilise le fait que 
$U_{z_{\phi^{-1}(r)}}^n\subseteq g_{w(r,q)}[V_{z_{\phi^{-1}(q)}}^{n-1}]$.\bigskip

\noindent 1.4.2. $\tilde A_\emptyset\subseteq T$.\bigskip

On a que $g_\emptyset^{-1}(g_{w(n,m)}[V_{z_{\phi^{-1}(m)}}^{n-1}])
\cap V_{z_{\phi^{-1}(n)}\lceil p}$ est non vide. On proc\`ede 
alors comme dans le cas 1.1.2 ; on a la condition (3) dans le cas 1.4.2.2, $r=n$, 
puisqu'on a $U_{z_{\phi^{-1}(n)}}^n\subseteq 
g_{w(r,q)}[V_{z_{\phi^{-1}(q)}}^{n-1}]$.\bigskip

\noindent\bf Cas 2.\rm\ $o=p+1$.\bigskip

\noindent 2.1. $z_{\phi^{-1}(m)}~{\cal R}~z_{\phi^{-1}(n)}$.\bigskip

 Soit $w\in \{\emptyset\}\cup\omega$ tel que 
$(N_{z_{\phi^{-1}(m)\lceil p}}\times 
N_{z_{\phi^{-1}(n)\lceil p}})\cap \mbox{Gr}(f_w)\not=\emptyset$. On peut supposer que 
$$\vert w\vert  = m(z_{\phi^{-1}(m)\lceil p},z_{\phi^{-1}(n)\lceil p}).$$ 
Par (iii), $\vert w(z_{\phi^{-1}(m)\lceil p},z_{\phi^{-1}(n)\lceil p})\vert  = \vert w\vert  = 0$ et on a  
$z_{\phi^{-1}(m)\lceil p} = z_{\phi^{-1}(n)\lceil p}$. Il y a deux cas.

\vfill\eject

\noindent 2.1.1. $\tilde A_\emptyset \subseteq Z$.\bigskip

 Le produit $U_{z_{\phi^{-1}(m)}}^{n-1}\times V_{z_{\phi^{-1}(m)}}^{n-1}$ rencontre 
$G(g_\emptyset )$, puisque $V_{z_{\phi^{-1}(m)}}^{n-1} = 
g_\emptyset[U_{z_{\phi^{-1}(m)}}^{n-1}]$. Comme ${G(g_\emptyset) = 
\overline{\bigcup_{n\in\omega} G(g_{n})}\setminus 
(\bigcup_{n\in\omega} G(g_{n}))}$, on peut trouver ${t\in \omega}$ minimal tel que 
$$(U_{z_{\phi^{-1}(m)}}^{n-1}\times V_{z_{\phi^{-1}(m)}}^{n-1})\cap G(g_{t})
\not=\emptyset\mbox{,}$$ 
et aussi on a ${\vert t\vert  = m(z_{\phi^{-1}(m)},z_{\phi^{-1}(n)}) = 1}$. On pose alors 
$\Phi(z_{\phi^{-1}(m)},z_{\phi^{-1}(n)}) := t$. Il y a 2 cas.\bigskip

\noindent 2.1.1.1. $\tilde A_{w(m,n)} \subseteq Z$.\bigskip

 On a alors $U_{z_{\phi^{-1}(n)\lceil p}} \cap g_\emptyset^{-1}(g_{w(m,n)}
[U_{z_{\phi^{-1}(m)}}^{n-1}\cap \tilde A_{w(m,n)}])\not=\emptyset$, et on raisonne 
comme en 1.1.1.\bigskip

\noindent 2.1.1.2. $\tilde A_{w(m,n)} \subseteq T$.\bigskip

 On a alors $U_{z_{\phi^{-1}(n)\lceil p}} \cap g_\emptyset^{-1}(g_{w(m,n)}^{-1}(
U_{z_{\phi^{-1}(m)}}^{n-1}))\not=\emptyset$, et on raisonne comme en 1.2.1.\bigskip

\noindent 2.1.2. $\tilde A_{\emptyset}\subseteq T$.\bigskip

 On a $U_{z_{\phi^{-1}(m)}}^{n-1} = g_\emptyset [V_{z_{\phi^{-1}(m)}}^{n-1}]$. On 
raisonne alors comme en 2.1.1 (on se r\'ef\`ere \`a 1.1.2 et 1.2.2 ; dans le cas analogue \`a 
1.1.2, on a $g_{w(m,n)}[U_{z_{\phi^{-1}(m)}}^{n-1}\cap \tilde A_{w(m,n)}]
\cap g_\emptyset^{-1}(U_{z_{\phi^{-1}(n)}\lceil p})\not=\emptyset$).\bigskip

\noindent 2.2. $z_{\phi^{-1}(n)}~{\cal R}~z_{\phi^{-1}(m)}$.\bigskip

 On raisonne comme en 2.1 : il y a deux cas.\bigskip
 
\noindent 2.2.1. $\tilde A_{\emptyset}\subseteq Z$.\bigskip

\noindent 2.2.1.1. $\tilde A_{w(n,m)} \subseteq Z$.\bigskip

 On a alors $U_{z_{\phi^{-1}(n)\lceil p}} \cap g_{w(n,m)}^{-1}(
V_{z_{\phi^{-1}(m)}}^{n-1})\not=\emptyset$, et on raisonne 
comme en 1.3.1.\bigskip

\noindent 2.2.1.2. $\tilde A_{w(n,m)} \subseteq T$.\bigskip

 On a alors 
$g_{w(n,m)}[V_{z_{\phi^{-1}(m)}}^{n-1}\cap \tilde A_{w(n,m)}]\not=\emptyset$, 
et on raisonne comme en 1.4.1.\bigskip

\noindent 2.2.2. $\tilde A_{\emptyset}\subseteq T$.\bigskip

 On raisonne comme en 2.1.2 et 2.2.1 (on se r\'ef\`ere \`a 1.3.2 et 1.4.2 ; dans le cas 
analogue \`a 1.4.2, on a $g_\emptyset^{-1}(g_{w(n,m)}[V_{z_{\phi^{-1}(m)}}^{n-1}\cap 
\tilde A_{w(n,m)}])\cap V_{z_{\phi^{-1}(n)}\lceil p}\not=\emptyset$).$\hfill\square$

\begin{thm} Il existe un bor\'elien $A_1$ de $2^\omega\times 2^\omega$ tel que 
pour tous espaces polonais $X$ et $Y$, et pour tout bor\'elien $A$ de $X\times Y$ qui est 
$\mbox{pot}(\borathree)$ et $\mbox{pot}(\bormthree)$, on a l'\'equivalence entre les conditions suivantes :\smallskip

\noindent (a) Le bor\'elien $A$ n'est pas $\mbox{pot}(\bormone)$.\smallskip

\noindent (b) Il existe des fonctions continues $u : 2^\omega \rightarrow X$ et 
$v : 2^\omega \rightarrow Y$ telles que $\overline{A_1} \cap (u\times v)^{-1}(A) = A_1$.\end{thm}

\vfill\eject

\noindent\bf D\'emonstration.\rm\ Soit $\Phi : 2^\omega \rightarrow Z_0$ un hom\'eomorphisme. 
On pose 
$$A_1 := (\Phi\times\Phi )^{-1}\left(\bigcup_{n>0} \mbox{Gr}(f_n)\right).$$
\noindent $\bullet$ Appliquons le th\'eor\`eme 2.3 \`a $X=Y=2^\omega$, $A=A_1$, $Z=T=Z_0$, 
$g=\mbox{Id}_{Z_0}$, $g_n = f_{n+1}$, et $u=v=\Phi^{-1}$. Ce th\'eor\`eme s'applique car\bigskip

\noindent - $A_1$ a ses coupes d\'enombrables, donc est $\mbox{pot}(\borathree)$ et $\mbox{pot}(\bormthree)$.\bigskip

\noindent - $(Z_0,(f_n)_{n>0})$ est une bonne situation (par 2.7), donc une situation g\'en\'erale 
puisque les domaines des $f_n$ sont non vides.\bigskip

\noindent - $\bigcup_{n>0} \mbox{Gr}(f_n) = (u\times v)^{-1}(A_1) = (u\times v)^{-1}(A_1)\cap
\overline{\bigcup_{n>0} \mbox{Gr}(f_n)}$.\bigskip 

 On a alors que $A_1$ n'est pas potentiellement ferm\'e, ce qui montre que (b) implique (a).\bigskip
 
\noindent $\bullet$ R\'eciproquement, si $A$ n'est pas potentiellement ferm\'e, le th\'eor\`eme 
2.3 nous fournit une situation g\'en\'erale $(Z,T,g_\emptyset,(g_n))$ et des injections 
continues $\tilde u : Z\rightarrow X$ et $\tilde v : T\rightarrow Y$ telles que 
$\overline{\bigcup_{n\in\omega} G(g_n)} \cap (\tilde u\times\tilde v)^{-1}(A) = 
\bigcup_{n\in\omega} G(g_n)$. Le th\'eor\`eme 2.8 nous fournit des fonctions continues 
$u':\omega^\omega\rightarrow Z$ et $v':\omega^\omega\rightarrow T$ telles que 
$$\overline{\bigcup_{n\in\omega} \mbox{Gr}(f_n)}\cap (u'\times v')^{-1}
\left(\bigcup_{n\in\omega} G(g_n)\right) = \bigcup_{n\in\omega} \mbox{Gr}(f_n).$$ 
Soit $\Psi: Z_0\rightarrow \omega^\omega$ l'injection canonique. On pose 
$$\left\{\!\!\begin{array}{ll}
u\!\!\!\!
& :=\tilde u\circ~u'\circ~\Psi\circ~\Phi\mbox{,}\cr
v\!\!\!\!
& :=\tilde v\circ~v'\circ~\Psi\circ~\Phi .
\end{array}
\right.$$
Alors $u$ et $v$ sont clairement continues, et on a clairement 
$A_1\subseteq (u\times v)^{-1}(A)$. Si $(\alpha,\beta)\in \overline{A_1}\setminus A_1$, 
$\alpha = \beta$ car $(\mbox{Gr}(f_n))_{n>0}$ converge vers ${\it\Delta} (Z_0)$, donc 
$$(u'(\Phi(\alpha )),v'(\Phi(\alpha )))\in\overline{\bigcup_{n\in\omega} G(g_n)}
\setminus\left(\bigcup_{n\in\omega} G(g_n)\right) \subseteq \check (\tilde u\times\tilde v)^{-1}(A).$$ 
Par cons\'equent, $(u(\alpha),v(\beta))\notin A$.$\hfill\square$\bigskip

\bf\noindent Remarques.\rm~(a) L'\'enonc\'e de ce th\'eor\`eme fournit un exemple de test $A_1$ 
dans $2^\omega \times 2^\omega$ ; mais la preuve nous donne aussi un test dans 
$Z_0\times Z_0$, et un autre dans $\omega^\omega \times \omega^\omega$. Nous donnerons une autre 
variante de cet exemple, qui se pr\^ete plus \`a g\'en\'eralisation, dans le paragraphe 
suivant.\bigskip

\noindent (b)  On peut d\'eterminer la complexit\'e d'un test comme dans l'\'enonc\'e 2.9 : $A_1$ est 
$D_2(\boraone)$ non $\mbox{pot}(\check D_2(\boraone ))$. En effet, avec l'ouvert  
$A = (2^\omega\times2^\omega )\setminus {\it\Delta} (2^\omega)$, qui n'est pas 
potentiellement ferm\'e (par la proposition 2.2 de [Le1]), on voit que $A_1$ est 
$D_2(\boraone)$. Mais un tel $A_1$ n'est pas $\mbox{pot}(\check D_2(\boraone ))$.

\vfill\eject

 En effet, avec $A = A_1$, on voit que $A_1$ n'est pas $\mbox{pot}(\bormone)$. Donc si $A_1$ \'etait 
$\mbox{pot}(\check D_2(\boraone ))$, on pourrait trouver une fonction continue 
$g:2^\omega\times2^\omega\rightarrow 2^\omega\times2^\omega$ telle que l'image 
r\'eciproque par $g$ de tout ensemble $\mbox{pot}(\boraone)$ soit $\mbox{pot}(\boraone)$, et telle que 
$g^{-1}(A_1)=(2^\omega\times2^\omega )\setminus {\it\Delta} (2^\omega)$ (cf [Le1], 
Cor. 4.14.(a)). On aurait donc, avec $l := (u\times v) \circ g$ et $A=D_0$, 
$g^{-1}(\overline{A_1})\cap l^{-1}(D_0) = 
(2^\omega\times2^\omega )\setminus {\it\Delta} (2^\omega)$ (cf [Le1], Ex. 3.6). Mais 
comme $g^{-1}(\overline{A_1})$ est ferm\'e de $2^\omega\times2^\omega$ et contient 
l'ouvert dense $(2^\omega\times2^\omega )\setminus {\it\Delta} (2^\omega)$, 
$g^{-1}(\overline{A_1}) = 2^\omega\times2^\omega$ et 
$l^{-1}(D_0)=(2^\omega\times2^\omega )\setminus {\it\Delta} (2^\omega)$, 
ce qui contredit la preuve du th\'eor\`eme 3.7 de [Le1].\bigskip

\bf\noindent (C) L'impossibilit\'e de l'injectivit\'e de la r\'eduction.\bigskip\rm

 Nous montrons maintenant qu'il n'est pas possible d'avoir $u$ et $v$ injectives dans 
le th\'eor\`eme 2.9. Cependant, il y a un cas o\`u on peut avoir l'injectivit\'e de la r\'eduction : 
quand $A = \bigcup_{n\in\omega} \mbox{Gr}(g_n)$, o\`u $(Z,(g_n))$ est une situation correcte (voir 
la d\'efinition ci-apr\`es).

\begin{defis} (1) On dit que $(Z,(g_n))$ est une $situation\ correcte$ si\smallskip

\noindent (a) $Z$ est un espace polonais parfait de dimension 0 non vide.\smallskip

\noindent (b) $g_n$ est un hom\'eomorphisme de domaine et d'image ouverts-ferm\'es de $Z$.\smallskip

\noindent (c) La suite $(\mbox{Gr}(g_n))$ converge vers la diagonale ${\it\Delta} (Z)$.\medskip

\noindent (2) On dit que $(Z,(g_n))$ est une $situation\ tr\grave es\ correcte$ si\smallskip

\noindent (a) $(Z,(g_n))$ est une situation correcte.\smallskip

\noindent (b) Pour toute suite finie d'entiers $p_0, ...,p_n$ et pour toute suite finie $\varepsilon _0, ...,\varepsilon _n$ d'\'el\'ements de $\{-1,1\}$, on a l'implication 
$$\exists~U\in \borone\lceil Z~~U\not=\emptyset~~\mbox{et}~~\forall~z\in U~~{g_{p_0}}^
{\varepsilon _0}...{g_{p_n}}^{\varepsilon _n}(z) = z~~\Rightarrow ~~\exists~i<n~~p_i = p_{i+1}~~
\mbox{et}~~\varepsilon _i = -\varepsilon_{i+1}.$$\end{defis}

 Par exemple, une bonne situation est une situation correcte.

\begin{lem} Soit $(Z,(g_n))$ une situation correcte. Alors il existe un 
ouvert-ferm\'e $C'_n$ du domaine de $g_n$ tel que si $g'_n := g_n\lceil C'_n$, 
$(Z,(g'_n))$ soit une situation tr\`es correcte.\end{lem}

\noindent\bf D\'emonstration.\rm\ Soit $(U_n)$ une base de la 
topologie de $Z$, et $V_n\subseteq U_n$ un ouvert-ferm\'e non vide de diam\`etre au plus 
$2^{-n}$. On construit par r\'ecurrence une suite injective $(q_k)$ d'entiers et un 
ouvert-ferm\'e non vide $C'_{q_k}$ de $C_{q_k}$ tels que\bigskip

\noindent (1) $\mbox{Gr}(g'_{q_k})\subseteq V_k^2$.\smallskip

\noindent (2) Pour toute suite finie $p_0, ...,p_n$ d'\'el\'ements de $\{q_0,...,q_k\}$ et pour toute 
suite finie $\varepsilon _0, ...,\varepsilon _n$ d'\'el\'e-ments de $\{-1,1\}$, on a l'implication 
$$\exists~U\in \borone\lceil Z~~U\not=\emptyset~\mbox{et}~ \forall~z\in U~~{g'_{p_0}}^
{\varepsilon _0}...{g'_{p_n}}^{\varepsilon _n}(z) = z~~\Rightarrow ~~\exists~i<n~~p_i = p_{i+1}~~
\mbox{et}~~\varepsilon _i = -\varepsilon_{i+1}.$$
(3) Il n'y a qu'un nombre fini de compositions d'\'el\'ements de $\{g'_{q_0},...,g'_{q_k},
{g'_{q_0}}^{-1},...,{g'_{q_k}}^{-1}\}$ ayant un domaine de d\'efinition non vide.\bigskip

\noindent $\bullet$ Admettons ceci r\'ealis\'e. Il restera \`a poser $C'_m := \emptyset$ si 
$m\notin \{q_k~/~k\in\omega\}$ pour avoir le lemme. Supposons donc la construction 
r\'ealis\'ee pour $l<k$.

\vfill\eject

 Enum\'erons l'ensemble fini $C$ des compositions d'\'el\'ements de 
$\{g'_{q_0},...,g'_{q_{k-1}},{g'_{q_0}}^{-1},...,{g'_{q_{k-1}}}^{-1}\}$ de domaine de d\'efinition non 
vide ne v\'erifiant pas la conclusion de l'implication de la condition (2) : $C = 
\{f_1,...,f_m\}$. On notera $D_{f_i}$ le domaine, n\'ecessairement ouvert-ferm\'e 
non vide, de $f_i$. Posons, pour $I\subseteq m$, 
$$O_I := \bigcap_{i\in I} D_{f_{i+1}} \cap \bigcap_{i\in m\setminus I} \check D_{f_{i+1}}.$$
Alors $(O_I)_{I\subseteq m}$ est une partition de $Z$ en ouverts-ferm\'es, et $\exists I\subseteq m$ tel que $V_k \cap O_I \not= \emptyset$. Posons 
$$O' := \{z\in V_k \cap O_I~/~\forall~i\in I~~f_{i+1}(z)\not= z\}.$$
Alors $O'$ est ouvert dense de $V_k \cap O_I$, donc par continuit\'e il existe un ouvert-ferm\'e 
non vide $O''$ de $O'$ tel que pour $i$ dans $I$ on ait $O''\cap f_{i+1}[O''] = \emptyset$. 
Choisissons $q_k\notin \{q_0,...,q_{k-1}\}$ tel que $\mbox{Gr}(g_{q_k})\cap (O''\times O'')
\not= \emptyset$, puis un ouvert-ferm\'e non vide $C'_{q_k}$ de $C_{q_k}\cap O''\cap 
g_{q_k}^{-1}(O'')$ tel que $C'_{q_k}\cap g_{q_k}[C'_{q_k}] = \emptyset$. \bigskip

\noindent $\bullet$ La condition (1) est clairement r\'ealis\'ee. Soient $p_0, ...,p_n$ une suite finie 
d'\'el\'ements de $\{q_0,...,q_k\}$, $\varepsilon _0, ...,\varepsilon _n$ une suite finie d'\'el\'ements de 
$\{-1,1\}$, $U$ un ouvert-ferm\'e non vide de $Z$ tel que pour $z$ dans $U$ on ait 
${g'_{p_0}}^{\varepsilon _0}...{g'_{p_n}}^{\varepsilon _n}(z) = z$, avec $p_i \not= p_{i+1}$ ou 
$\varepsilon _i = \varepsilon_{i+1}$ pour tout $i<n$. Alors il existe $i\leq n$ tel que $p_i = q_k$, 
par hypoth\`ese de r\'ecurrence ; montrons qu'un tel $i$ est unique. Si tel n'est pas 
le cas, dans la composition appara\^\i t ${g'_{q_k}}^{\varepsilon _j}{g'_{p_{j+1}}}^{\varepsilon _{j+1}}
...{g'_{p_s}}^{\varepsilon _s}{g'_{q_k}}^{\varepsilon _{s+1}}$, avec $p_l\not= q_k$ si $j<l\leq s$. 
On a ${g'_{q_k}}^{\varepsilon _{s+1}}...{g'_{p_n}}^{\varepsilon _n}(z)\in O''\subseteq O_I$. Soit 
$r < m$ tel que 
$f_{r+1} = {g'_{p_{j+1}}}^{\varepsilon _{j+1}}...{g'_{p_s}}^{\varepsilon _s}$ ; alors 
$r\in I$ car $f_{r+1}({g'_{q_k}}^{\varepsilon _{s+1}}...{g'_{p_n}}^{\varepsilon _n}(z))$ est d\'efini et 
${g'_{q_k}}^{\varepsilon _{s+1}}...{g'_{p_n}}^{\varepsilon _n}(z)\in O_I$. Donc 
$f_{r+1}({g'_{q_k}}^{\varepsilon _{s+1}}...{g'_{p_n}}^{\varepsilon _n}(z))\notin O''$, ce qui est absurde. 
D'o\`u l'unicit\'e de $i$. On en d\'eduit l'existence d'un ouvert-ferm\'e non vide de $Z$ 
sur lequel $g'_{q_k}$ co\"\i ncide avec l'une des fonctions $f_t$ ou avec l'identit\'e. Mais 
ceci contredit le choix de $C'_{q_k}$ et de $O''$. D'o\`u la condition (2).\bigskip

 Posons $H := \{g'_{q_k},{g'_{q_k}}^{-1},\mbox{Id}_{C'_{q_k}},\mbox{Id}_{g'_{q_k}[C'_{q_k}]}\}$. Alors 
les seules compositions d'\'el\'ements de $H$ ayant un domaine d\'efinition non vide sont les 
\'el\'ements de $H$. Par suite, les seules compositions d'\'el\'ements de $\{g'_{q_0},...,
g'_{q_k},{g'_{q_0}}^{-1},...,{g'_{q_k}}^{-1}\}$ ayant un domaine de d\'efinition non vide sont de la 
forme $h$, $f$, $hf$, $fh$, ou $fhg$, o\`u $f$ et $g$ sont des compositions d'\'el\'ements de 
$$\{g'_{q_0},...,g'_{q_{k-1}},{g'_{q_0}}^{-1},...,{g'_{q_{k-1}}}^{-1}\}$$ 
de domaine de d\'efinition non vide (on raisonne comme pr\'ec\'edemment pour voir qu'il n'y a pas plus d'un \'el\'ement de $H$ dans une composition) ; elles sont donc en nombre fini. D'o\`u la condition (3) et le 
lemme.$\hfill\square$

\begin{thm} Soit $(Z,(g_n))$ une situation correcte. Alors il existe $u : Z_0\rightarrow Z$ injective continue telle que 
$$\overline{\bigcup_{n>0} \mbox{Gr}(f_n)}\cap (u\times u)^{-1}
\left(\bigcup_{n\in\omega} \mbox{Gr}(g_n)\right) = \bigcup_{n>0} \mbox{Gr}(f_n).$$\end{thm}

\noindent\bf D\'emonstration.\rm\ On utilisera des notations analogues \`a celles de la preuve du th\'eor\`eme 2.6, et le m\^eme sch\'ema de d\'emonstration. Les nuances sont les suivantes. Soit $(g'_n)$ fournie par le lemme 2.11.

\vfill\eject

\noindent $\bullet$ On va construire\bigskip

\noindent - Une suite $(U_s)_{s\in \bigcup_{p\in\omega} \Pi_{n<p}~{\cal A}_n}$ d'ouverts-ferm\'es non 
vides de $Z$.\smallskip

\noindent - Une fonction $\Phi :\bigcup_{p\in\omega} [(\Pi_{n<p}~{\cal A}_n)\times (\Pi_{n<p}~
{\cal A}_n )] \rightarrow \{\emptyset\}\cup \omega$.\bigskip

 On demande \`a ces objets de v\'erifier 
$$\begin{array}{ll}
& (i)~~~U_{s^\frown i}\subseteq U_s\cr 
& (ii)~~{\delta} (U_{s^\frown i})\leq 2^{-\vert s\vert -1} \cr 
& (iii)~s~{\cal R}~t \Rightarrow \left\{\!\!\!\!\!\!
\begin{array}{ll}
& \vert w(s,t)\vert  = m(s,t) \cr 
& U_t = g'_{w(s,t)} [U_s] 
\end{array}
\right.\cr
& (iv)~~U_{s^\frown n}\cap U_{s^\frown m} = \emptyset ~\mbox{si}~n\not= m
\end{array}$$
On peut d\'efinir $u:Z_0\rightarrow Z$ comme dans la preuve du th\'eor\`eme 2.6, et $u$ est injective 
continue. Comme $\mbox{Gr}(g'_n)\subseteq \mbox{Gr}(g_n)\subseteq \check {\it\Delta} (Z)$, la construction 
permet d'avoir le th\'eor\`eme.\bigskip

\noindent$\bullet$ Montrons donc que la construction est possible. Admettons avoir 
construit $U_s$ et $\Phi (s,t)$ pour $\vert s\vert $, $\vert t\vert \leq p$ 
v\'erifiant (i)-(iv), et soient $s\in \Pi_{n<p}~{\cal A}_n$ et $i\in {\cal A}_p$. On va 
construire par r\'ecurrence sur $n\in\{1,...,p_0+...+p_q\}$, et 
pour $k\in \{1,...,n\}$, des ouverts-ferm\'es non vides $U^n_{z_{\phi^{-1}(k)}}$ de $Z$. 
On demande aux ouverts-ferm\'es de v\'erifier 
$$\begin{array}{ll} 
& (1)~U^n_{z_{\phi^{-1}(k)}} \subseteq U_{z_{\phi^{-1}(k)}\lceil p}\cr 
& (2)~{\delta} (U^n_{z_{\phi^{-1}(k)}})\leq 2^{-p-1}\cr
& (3)~\mbox{Si}~k,l\in\{1,...,n\}~\mbox{et}~z_{\phi^{-1}(k)}~{\cal R}~z_{\phi^{-1}(l)}
\mbox{,~alors} \cr 
& ~~~~~~-\vert w(k,l)\vert  = m(z_{\phi^{-1}(k)},z_{\phi^{-1}(l)})  \cr
& ~~~~~~-U_{z_{\phi^{-1}(l)}}^n = g'_{w(k,l)}[U_{z_{\phi^{-1}(k)}}^n] \cr 
& (4)~U^{n+1}_{z_{\phi^{-1}(k)}} \subseteq 
U^n_{z_{\phi^{-1}(k)}}~\mbox{si}~k\in\{1,...,n\}\cr
& (5)~U^n_{z_{\phi^{-1}(k)}}\cap U^n_{z_{\phi^{-1}(l)}}=\emptyset~\mbox{si}~k,l\in
\{1,...,n\}~\mbox{et}~k\not= l
\end{array}$$
 On voit comme dans la preuve du th\'eor\`eme 2.6 que cette construction est suffisante 
((iv) r\'esulte du fait que $\Pi_{n\leq p}~{\cal A}_n$ est une classe pour $\cal E$, 
comme on l'a vu en 2.7).\bigskip

 Pour avoir la condition (5), on utilise le fait, vu en 2.7, que $(Z_0,(f_n)_{n>0})$ est 
une tr\`es bonne situation. On commence par assurer les conditions (1)-(4) comme dans la 
preuve du th\'eor\`eme 2.6. Soient donc $s$ et $t$ dans $\Pi_{n\leq p}~{\cal A}_n$, avec 
$s\not= t$ ; comme ce produit est une $\cal E$-classe, il existe une $\cal T$-cha\^\i ne 
$c$ de longueur minimale telle que $c(0)=s$ et $c(\vert c\vert -1)=t$. Soient $n := \vert c\vert -2$, 
$p_i$ des entiers et $\varepsilon_i$ dans $\{-1,1\}$ tels que pour tout $\alpha$ dans 
$\Pi_{n>p}~{\cal A}_n$, ${f_{p_i}}^{\varepsilon_i}(c(n-i)^\frown\alpha )=c(n-i+1)^\frown\alpha$ 
(ces objets existent par minimalit\'e de $n$). Par minimalit\'e de $n$ encore, on a 
$p_i\not=p_{i+1}$ ou $\varepsilon_i = \varepsilon_{i+1}$ pour tout $i<n$. Comme $(Z_0,(g'_n))$ est 
une situation tr\`es correcte, on peut trouver $x\in U_s$ tel que 
${g'_{p_0}}^{\varepsilon_0}...{g'_{p_n}}^{\varepsilon_n}(x) \not= x$ ; par suite, il existe un 
voisinage ouvert-ferm\'e $U'_s$ de $x$, inclus dans $U_s$, tel que 
${g'_{p_0}}^{\varepsilon_0}...{g'_{p_n}}^{\varepsilon_n}[U'_s]\cap U'_s=\emptyset$. on construit 
alors \`a nouveau des ouverts-ferm\'es $U'_r$, pour $r$ dans $\Pi_{n\leq p}~{\cal A}_n$, 
par r\'ecurrence sur $d(r,s)$, comme dans la preuve du th\'eor\`eme 2.6. On a alors $U'_r
\subseteq U_r$ et l'unicit\'e de la $\cal T$-cha\^\i ne allant d'une suite \`a l'autre (qui  
r\'esulte de la condition (b) de la d\'efinition d'une tr\`es bonne situation) montre que 
$U'_s\cap U'_t=\emptyset$. En un nombre fini d'\'etapes, on obtient donc la condition (5) 
en plus des autres conditions (1)-(4).$\hfill\square$

\vfill\eject

\begin{lem} Soient $H$, $K$ des bor\'eliens de $Z_0$, et $B$ un bor\'elien de 
$\bigcup_{n>0}~\mbox{Gr}(f_n)$ inclus dans $H\times K$ et non $\mbox{pot}(\bormone)$. Alors il existe un bor\'elien $Z$ de $H\cap K$, une topologie $\tau$ sur $Z$, et un bor\'elien $C_n$ de $D_{f_n}\cap Z$ tels que\smallskip

\noindent (a) La topologie $\tau$ est plus fine que la topologie initiale de $Z$.\smallskip

\noindent (b) Le graphe de la restriction de $f_n$ \`a $C_n$ est inclus dans $B$. \smallskip

\noindent (c) Le couple $((Z,\tau),(f_n\lceil C_n)_{n>0})$ est une situation correcte.\end{lem}

\noindent\bf D\'emonstration.\rm\ Comme $B\subseteq \bigcup_{n>0}~\mbox{Gr}(f_n)$, 
on peut \'ecrire que $B = \bigcup_{n>0}~\mbox{Gr}(f_n\lceil E_n)$, o\`u $E_n$ est bor\'elien de $Z_0$. On peut supposer, pour simplifier l'\'ecriture, que $Z_0$ est r\'ecursivement pr\'esent\'e et que $H$, $K$, $B$, 
$E_n$ et $f_n\lceil E_n$ sont $\Borel$. On notera 
${\it{\it\Sigma}}$ (respectivement ${\it\Delta}$) la topologie engendr\'ee par les $\Ana$ (respectivement
 $\Borel$) de $Z_0$. Posons 
$$\begin{array}{ll} 
\Omega\!\!\!\! & := \{x\in Z_0~/~\omega_1^x = \omega_1^{\mbox{CK}}\}\mbox{,}\cr 
D\!\!\!\! & := \{x\in Z_0~/~x\notin\Borel\}\mbox{,}\cr
Z\!\! \!\!& := \{x\in H\cap K\cap D\cap\Omega~/~(x,x)\in \overline{B}^{{\it\Sigma}\times{\it\Sigma}}\}.
\end{array}$$
Alors $\Omega$ muni de ${\it\Sigma}$ est polonais, donc $\Omega$ est bor\'elien de $Z_0$ et 
$Z$ aussi (puisque $B$ est $\Ana$, une double application du th\'eor\`eme de s\'eparation nous 
donne que $\overline{B}^{{\it\Sigma}\times{\it\Sigma}} = 
\overline{B}^{{\it\Delta}\times{\it\Delta}}$, et on utilise le fait que ${\it\Delta}$ soit polonaise 
pour voir que $\overline{B}^{{\it\Sigma}\times{\it\Sigma}}$ est bor\'elien). On d\'efinit $\tau$ comme 
\'etant la restriction de ${\it\Sigma}$ \`a $Z$. Comme $Z$ est $\Ana$, c'est un ouvert de $\Omega$ 
muni de ${\it\Sigma}$, donc $(Z,\tau)$ est polonais, de dimension $0$ car traces des $\Ana$ 
sur $\Omega$ sont ouverts-ferm\'es de $\Omega$ muni de ${\it\Sigma}$, parfait car inclus dans 
l'ensemble co-d\'enombrable $D$ (je renvoie le lecteur \`a [Lo3] pour les preuves 
des propri\'et\'es de $\Omega$, ${\it\Delta}$ et ${\it\Sigma}$ utilis\'ees ici et non d\'emontr\'ees).\bigskip

 On va montrer que $Z\not= \emptyset$. Le bor\'elien $B$ n'est pas $\mbox{pot}(\bormone)$, donc 
$B\cap D^2$ non plus puisqu'un bor\'elien de projection d\'enombrable est $\mbox{pot}(\borone)$ 
(cf [Le1], remarque 2.1). Soit $\alpha$ tel que $D$ soit $\Borel (\alpha )$. 
Alors $D^2\cap \check B \cap \overline{B\cap D^2}^{{\it\Delta} (\alpha) \times{\it\Delta} 
(\alpha )} \not= \emptyset$, donc contient $(x,y)$. On peut donc trouver $(x_n,y_n)$ 
dans $\mbox{Gr}(f_{p_n}\lceil E_{p_n})$, o\`u $p_n\in\omega$, tel que $(x_n,y_n)$ converge vers 
$(x,y)$, pour ${\it\Delta} (\alpha )\times{\it\Delta} (\alpha )$. Or le graphe de $f_p\lceil E_p$ 
est ferm\'e dans $Z_0$ muni de ${\it\Delta} (\alpha )\times{\it\Delta} (\alpha )$ ; on 
peut donc supposer que la suite $(p_n)$ cro\^\i t strictement vers l'infini, car 
$(x,y)\notin B$. Par suite $y=x$, et comme $H$ et $K$ sont ferm\'es pour 
${\it\Delta} (\alpha )$, on a $x\in H\cap K$. On en d\'eduit que 
$(H\cap K \cap D)^2\cap \check B \cap \overline{B\cap D^2}^
{{\it\Delta}\times{\it\Delta}} \not= \emptyset$. Mais ce $\Ana$ non vide, qui vaut 
$(H\cap K \cap D)^2\cap \check B \cap \overline{B\cap D^2}^{{\it\Sigma}\times{\it\Sigma}}$, 
rencontre $\{(x',y')\in Z_0\times Z_0~/~\omega_1^{(x',y')} = \omega_1^{\mbox{CK}}\}\subseteq 
\Omega^2$, en $(x',y')$, et en raisonnant comme pr\'ec\'edemment, on voit que $x'=y'\in Z$.
\bigskip

 Soit $(z_n)$ une suite dense de $(Z,\tau)$, et $Q_n$ un $\Ana$ non vide contenant 
$z_n$ inclus dans la boule ouverte de centre $z_n$ et de rayon $2^{-n}$, au sens de 
$(Z,\tau)$. Alors on construit par r\'ecurrence une suite injective d'entiers 
$(q_n)$ telle que $Q_n^2\cap \mbox{Gr}(f_{q_n}\lceil E_{q_n})\not=\emptyset$. On pose alors 
$C_{q_n} := E_{q_n}\cap Q_n\cap f_{q_n}^{-1}(Q_n)$, et $C_m := \emptyset$ si $m\notin 
\{q_n~/~n\in\omega\}$. Les ensembles $C_n$ sont ouverts-ferm\'es de $(Z,\tau)$, donc 
bor\'eliens de $Z_0$ muni de sa topologie initiale. Les ensembles $f_n[C_n]$ sont $\Ana$ 
donc ouverts-ferm\'es de $(Z,\tau)$. Les restrictions de $f_n$ \`a $C_n$ sont clairement 
des hom\'eomorphismes car tous les objets intervenant sont $\Ana$. Leurs graphes sont 
dans $B$ car $C_n\subseteq E_n$. Enfin, la suite $(\mbox{Gr}(f_n\lceil C_n))_{n>0}$ converge vers 
${\it\Delta} (Z)$ par construction des $C_n$.$\hfill\square$
 
\begin{lem} Soit $(h_n)_{n>0}$ une suite de fonctions continues et ouvertes de 
domaine et d'image $2^\omega$ telle que ${\it\Delta} (2^\omega) \subseteq \overline{
\bigcup_{n>0} \mbox{Gr}(h_n)}$ et $h_n(\alpha )\not= \alpha$ pour tout $\alpha$ d'un ouvert 
dense de $2^\omega$. Alors il existe un ouvert-ferm\'e $C_n$ de $2^\omega$ tel que si 
$h'_n := h_n\lceil C_n$, on ait\smallskip

\noindent (a) $(2^\omega,2^\omega,Id_{2^\omega},(h'_n)_{n>0})$ est une situation g\'en\'erale (\`a ceci 
pr\`es que certains $C_n$ peuvent \^etre vides).\smallskip

\noindent (b) Pour tout $n>0$, on a $C_n\cap h'_n [C_n] = \emptyset$.\end{lem}

\noindent\bf D\'emonstration.\rm\ Soit $(\alpha_n)$ une suite dense de $2^\omega$. On 
va construire, en plus des $C_n$, une suite 
d'entiers strictement croissante $(m(n))_n$. On choisit $x$ et $y$ dans 
${\cal B}(\alpha_n,2^{-n-1}[$ tels que $y\not= x$ et $y = h_{m(n)}(x)$, o\`u $m(n)$ est un 
entier strictement sup\'erieur \`a $\mbox{max}_{p<n}~m(p)$. On choisit un voisinage 
ouvert-ferm\'e $\tilde C_{m(n)}$ de $x$ tel que
$${(\tilde C_{m(n)}\times h_{m(n)}[\tilde C_{m(n)}])\cap {\it\Delta} (Z)=\emptyset}\mbox{,}$$  
$${\tilde C_{m(n)}\subseteq 
{\cal B}(\alpha_n,2^{-n-1}[~\cap ~h_{m(n)}^{-1}({\cal B}(\alpha_n,2^{-n-1}[)}.$$ 
On choisit un ouvert-ferm\'e non 
vide $D_{m(n)}$ de $h_{m(n)}[\tilde C_{m(n)}]$, et puis enfin on pose 
$${C_{m(n)} := \tilde C_{m(n)}\cap h_{m(n)}^{-1}(D_{m(n)})}.$$ 
Il reste \`a poser $C_m := \emptyset$ si ${m\notin\{m(n)~/~n\in\omega\}}$.$\hfill\square$\bigskip

\noindent\bf Notation. \rm~On pose, pour $n>0$,
$$h_n : \left\{\!\!
\begin{array}{ll}
2^\omega\!\!\!\! 
& \rightarrow 2^\omega \cr 
\alpha 
& \mapsto \left\{\!\!
\begin{array}{ll}
\omega\!\!\!\! 
& \rightarrow 2 \cr 
k 
& \mapsto \left\{\!\!\!\!\!\!
\begin{array}{ll}
& \alpha (k)~\mbox{si}~k\not\equiv -1~(2^n)\mbox{,}\cr\cr 
& \alpha (2^{n+1}q-2^n-1)~\mbox{si}~k=2^nq-1.
\end{array}
\right.
\end{array}
\right.
\end{array}
\right.$$

\begin{lem} Les fonctions $h_n$ v\'erifient les conditions suivantes :\medskip

\noindent - Ce sont des surjections continues et ouvertes.\smallskip

\noindent - ${\it\Delta} (2^\omega) \subseteq \overline{\bigcup_{n>0} \mbox{Gr}(h_n)}$.\smallskip

\noindent - Pour tout $n>0$, $h_n (\alpha) \not= \alpha$ sur un ouvert dense de $2^\omega$.\smallskip

\noindent - $0<p<n~\Rightarrow~\forall\alpha\in2^\omega~~h_p(h_n(\alpha )) = h_p (\alpha )$.\end{lem}

\noindent\bf D\'emonstration.\rm\ Les deux premi\`eres conditions sont 
clairement r\'ealis\'ees car $(h_n)_{n>0}$ converge uniform\'ement vers $\mbox{Id}_{2^\omega}$. Si 
$s\in \omega^{<\omega}$ et $n>0$, on peut trouver un entier $q$ tel que 
$2^nq-1\geq\vert s\vert $ et $2^{n+1}q-2^n-1\geq\vert s\vert $ ; par cons\'equent, on trouve $\alpha$ 
dans $2^\omega$ tel que $h_n (s^\frown \alpha)\not= s^\frown \alpha$. 
D'o\`u la troisi\`eme condition. Il reste \`a voir la quatri\`eme :\bigskip

\noindent\bf Cas 1.\rm\ $k \not\equiv -1~(2^p)$.\bigskip

 On a alors $h_p(h_n(\alpha ))(k) = h_n(\alpha )(k) = \alpha (k)$, car 
$k \not\equiv -1~(2^n)$ sinon on trouve $q$ tel que $k =2^nq-1 = 2^p(2^{n-p}q)-1$. Par 
ailleurs, $h_p(\alpha )(k) = \alpha (k)$.\bigskip

\noindent\bf Cas 2.\rm\ $k \equiv -1~(2^p)$.\bigskip

 On a alors $h_p(h_n(\alpha ))(k) = h_n(\alpha )(2^{p+1}q-2^p-1) = \alpha 
(2^{p+1}q-2^p-1)$ car on a que ${2^{p+1}q-2^p-1\not\equiv -1~(2^n)}$ sinon on trouve 
$q'$ tel que $2^{p+1}q-2^p-1 = 2^nq'-1$, ce qui entra\^\i ne que $2q-1 = 2^{n-p}q'$. Par ailleurs, 
$h_p(\alpha )(k) = \alpha (2^{p+1}q-2^p-1)$.$\hfill\square$

\vfill\eject

\begin{thm} Il n'est pas possible d'avoir $u$ et $v$ injectives dans le th\'eor\`eme 2.9.\end{thm}

\noindent\bf D\'emonstration.\rm\ Raisonnons par l'absurde : il existe un test $B_1$ avec 
injectivit\'e des fonctions de r\'eduction. Avec $A = (2^\omega\times2^\omega )\setminus 
{\it\Delta} (2^\omega )$, on voit que $B_1$ est $D_2(\boraone)$. Avec $A= B_1$, on voit que 
$B_1\notin\mbox{pot}(\bormone )$. Avec $A = \bigcup_{n>0}~\mbox{Gr}(f_n)$, on trouve des 
injections continues $u$ et $v$ de $2^\omega$ dans $Z_0$ telles que $\overline{B_1}\cap 
(u\times v)^{-1}(\bigcup_{n>0}~\mbox{Gr}(f_n))= B_1$. Posons $B := (u\times v)[B_1]$, $H := u[2^\omega]$, 
$K := v[2^\omega]$. Alors $B\notin\mbox{pot}(\bormone )$, sinon $B_1$ le serait. 
Le lemme 2.13 peut donc s'appliquer et nous fournit, avec le th\'eor\`eme 2.12, une 
injection continue $\tilde u:Z_0 \rightarrow (Z,\tau)$ telle que 
$$\overline{\bigcup_{n>0} \mbox{Gr}(f_n)}\cap (\tilde u\times \tilde u)^{-1}
\left(\bigcup_{n>0} \mbox{Gr}(f_n\lceil C_n)\right) = \bigcup_{n>0} \mbox{Gr}(f_n).$$
Autrement dit, il existe une injection continue $u':Z_0\rightarrow u[2^\omega]\cap v[2^\omega]$ telle que 
$$\overline{\bigcup_{n>0}~\mbox{Gr}(f_n)}\cap (u'\times u')^{-1}(B) = \bigcup_{n>0}~\mbox{Gr}(f_n).$$ 
Posons $Z' := u'[Z_0]$, puis 
$${\cal U}: \left\{\!\!
\begin{array}{ll}
Z_0\!\!\!\! 
& \rightarrow 2^\omega \cr 
z 
& \mapsto u^{-1}(u'(z))
\end{array}
\right.\mbox{,}~~{\cal V}: \left\{\!\!
\begin{array}{ll}
Z_0\!\!\!\!
& \rightarrow 2^\omega \cr 
z 
& \mapsto v^{-1}(u'(z))
\end{array}.
\right.$$ 
 Alors ${\cal U}$ et ${\cal V}$ sont des injections continues telles que 
$$\overline{\bigcup_{n>0}~\mbox{Gr}(f_n)}\cap ({\cal U}\times {\cal V})^{-1}(B_1)=
\bigcup_{n>0}~\mbox{Gr}(f_n).$$
En effet, si 
$(\alpha,\beta)\in \overline{\bigcup_{n>0}~\mbox{Gr}(f_n)}\cap 
({\cal U}\times {\cal V})^{-1}(B_1)\setminus 
(\bigcup_{n>0}~\mbox{Gr}(f_n))$, $\alpha = \beta$. On en d\'eduit, puisque 
$(u({\cal U}(\alpha)),v({\cal V}(\alpha)))\notin \bigcup_{n>0}~\mbox{Gr}(f_n)$, que 
$({\cal U}(\alpha),{\cal V}(\beta))\notin B_1$, ce qui est absurde.\bigskip

 On en d\'eduit que si $A$ est $\mbox{pot}(\borathree)$ et $\mbox{pot}(\bormthree)$, mais pas 
potentiellement ferm\'e de $X\times Y$, on peut trouver des injections continues 
$U : Z_0 \rightarrow X$ et $V : Z_0 \rightarrow Y$ telles que 
$$\overline{\bigcup_{n>0}~\mbox{Gr}(f_n)}\cap (U\times V)^{-1}(A)=\bigcup_{n>0}~\mbox{Gr}(f_n).$$ 
Il suffit en effet de composer les 
injections continues qui seraient fournies par une version injective du th\'eor\`eme 2.9 
avec ${\cal U}$ et ${\cal V}$.\bigskip

 Par le lemme 2.15, on peut appliquer le lemme 2.14 \`a $Z=2^\omega$ et \`a $(h_n)_{n>0}$, 
ce qui fournit une situation g\'en\'erale ${(2^\omega,2^\omega,\mbox{Id}_{2^\omega},(h'_n)_{n>0})}$. 
On a, par application du th\'eor\`eme 2.3, que $\bigcup_{n>0}~\mbox{Gr}(h'_n)$ n'est pas 
potentiellement ferm\'e. On peut donc trouver des injections continues $U$ et $V$ de 
$Z_0$ dans $2^\omega$ telles que 
$\overline{\bigcup_{n>0}~\mbox{Gr}(f_n)}\cap (U\times V)^{-1}(\bigcup_{n>0}~\mbox{Gr}(h'_n))=
\bigcup_{n>0}~\mbox{Gr}(f_n)$. Montrons que ceci 
est impossible, ce qui constituera la contradiction cherch\'ee.\bigskip

 On raisonne par l'absurde. On commence par remarquer que $U=V$, puisque si $z\in Z_0$, alors 
$(z,z)\in\overline{\bigcup_{n>0}~\mbox{Gr}(f_n)}\setminus (\bigcup_{n>0}~\mbox{Gr}(f_n))$, donc 
$$(U(z),V(z))\in\overline{\bigcup_{n>0}\mbox{Gr}(h'_n)}\setminus\left(\bigcup_{n>0} \mbox{Gr}(h'_n)\right) = {\it\Delta} (2^\omega ).$$ 

\vfill\eject

Posons, pour $n>0$ et $m>0$, 
$$A_m^n := \{\alpha\in D_{f_n}~/~U(\alpha )\in C'_m~\mbox{et}~U(f_n(\alpha ))\in h'_m[C'_m]\}.$$
Alors $D_{f_n}\subseteq\bigcup_{m>0} A^n_m$, donc on peut trouver $m$ tel que 
$A^1_m\not=\emptyset$. On a, pour $\alpha\in A^1_m$, 
$$(U(\alpha ),U(f_1(\alpha )))\in \mbox{Gr}(h'_m)\mbox{,}$$ 
donc $h'_m(U(\alpha )) = U(f_1(\alpha ))$. Par cons\'equent, la restriction de $h'_m$ \`a 
$U[A^1_m]$ est injective. \bigskip

 Soit $\alpha_0$ dans $A^1_m$ ayant une infinit\'e de $1$ ($\alpha_0$ existe puisque 
$A^1_m$ est ouvert de $Z_0$). Notons $(n_k)_{k\in\omega}$ la suite des entiers tels 
que $\alpha_0\in D_{f_{n_k}}$. On peut trouver $p_k>0$ tel que 
$$U(f_{n_k}(\alpha_0)) = h'_{p_{n_k}}(U(\alpha_0)).$$ 
Alors pour $k$ assez grand, on a $f_{n_k}(\alpha_0)\in A^1_m$. Il existe $k_0$ assez grand tel que $p_{n_{k_0}}>m$, sinon on trouve $r\leq m$ tel que $U(f_{n_k}(\alpha_0)) = h'_r(U(\alpha_0))$ pour une infinit\'e de $k$. Comme $(f_{n_k}(\alpha_0))_{k\in\omega}$ converge vers $\alpha_0$, on a alors que $U(\alpha_0) = h'_r(U(\alpha_0))\in C_r\cap h'_r[C_r]=\emptyset$, ce qui est absurde. 
On a donc $h'_m(h'_{p_{n_{k_0}}}(U(\alpha_0))) = h'_m(U(\alpha_0))$. Comme 
$U(\alpha_0)$ et $h'_{p_{n_{k_0}}}(U(\alpha_0))=U(f_{n_{k_0}}(\alpha_0))$ sont dans 
$U[A^1_m]$, on a par injectivit\'e que $U(\alpha_0) = h'_{p_{n_{k_0}}}(U(\alpha_0))\in 
C_{p_{n_{k_0}}}\cap h'_{p_{n_{k_0}}}[C_{p_{n_{k_0}}}]=\emptyset$, ce qui est absurde.
$\hfill\square$\bigskip

\noindent\bf
{\Large 3 Un test pour les ensembles non potentiellement diff\'erence transfinie 
d'ouverts.}\bigskip\rm

 Dans le paragraphe 2, nous avons vu le r\^ole important jou\'e par l'arbre $\{\emptyset\}\cup\omega$. 
Nous g\'en\'eralisons maintenant cette notion d'arbre.\bigskip

\noindent\bf Notations. \rm~Dans toute la suite, $\xi$ d\'esignera un ordinal d\'enombrable.\bigskip

\noindent $\bullet$ On d\'efinit $\Psi_\xi : \omega^{<\omega}\rightarrow 
\{-1\}\cup (\xi+1)$ par r\'ecurrence sur $\vert s\vert $, o\`u $s\in \omega^{<\omega}$ : 
$\Psi_\xi (\emptyset)\! =\! \xi$ et 
$$\Psi_\xi(s^\frown n)\! =\! \left\{\!\!\!\!\!\!
\begin{array}{ll}
& \bullet -1~\mbox{si}~\Psi_\xi(s)\leq 0\mbox{,}\cr 
& \bullet\ \theta ~\mbox{si}~\Psi_\xi(s)=\theta+1\mbox{,}\cr 
& \bullet~\mbox{un~ordinal~impair~de}~\Psi_\xi(s)~\mbox{tel~que~la~suite}~
(\Psi_\xi(s^\frown n))_n~\mbox{soit~cofinale}\cr 
& \mbox{dans}~\Psi_\xi(s)~\mbox{et~strictement~croissante~si}~\Psi_\xi(s)~
\mbox{est~limite~non~nul.}
\end{array}
\right.$$
On d\'efinit alors des arbres : $T_\xi\! :=\! \{s\in\omega^{<\omega}~/~\Psi_\xi(s)\not= -1\}$ 
et ${T'_\xi\! :=\! \{s\in T_\xi~/~\Psi_\xi(s)\not= 0\}}$. Ces arbres sont bien s\^ur bien fond\'es, 
et la hauteur de $T_\xi$ (resp. $T'_\xi$) est $1+\xi$ (resp. $\xi$).\bigskip

\noindent $\bullet$ Soient $(A_s)_{s\in T_\xi}$, $(B_s)_{s\in T_\xi}$, 
$Z$ et $T$ des ensembles, 
avec $A_s\times B_s \subseteq Z\times T$ ou ${A_s\times B_s \subseteq T\times Z}$, et 
$f_s\! :\! A_s\!\rightarrow\! B_s$ des fonctions. On note 
$B_p := \bigcup_{s\in T_\xi,~\vert s\vert ~\mbox{paire}}~ G(f_s)$ et 
$B_i := \bigcup_{s\in T_\xi,~\vert s\vert ~\mbox{impaire}}~ G(f_s)$. 

\vfill\eject

\begin{defis} (1) On dit que $(Z,T,(f_s)_{s\in T_\xi})$ est une $situation\ g\acute en\acute erale$ si \smallskip

\noindent (a) $Z$ et $T$ sont des espaces polonais parfaits de dimension 0.\smallskip

\noindent (b) Les $f_s$ sont des fonctions continues et ouvertes de domaine ouvert-ferm\'e 
non vide de $Z$ et d'image ouverte-ferm\'ee de $T$, ou de domaine ouvert-ferm\'e non vide 
de $T$ et d'image ouverte-ferm\'ee de $Z$.\smallskip
 
\noindent (c) La suite $(G(f_{s^\frown n}))_n$ converge vers $G(f_s)$ si $s\in T'_\xi$.\smallskip

\noindent (d) $\overline{B_p} = B_p\cup_{\mbox{disj.}} B_i$ si $\xi$ est pair, et $\overline{B_i} = 
B_p\cup_{\mbox{disj.}} B_i$ si $\xi$ est impair.\medskip

\noindent (2) On dit que $(Z,(f_s)_{s\in T_\xi})$ est une $bonne\ situation$ si\smallskip

\noindent (a) $Z$ est un ferm\'e parfait non vide de $\omega^\omega$.\smallskip

\noindent (b) $f_s$ est un hom\'eomorphisme de domaine et d'image ouverts-ferm\'es non vides de $Z$, 
et les graphes des $f_s$ sont deux \`a deux disjoints. De plus, 
$\alpha\leq_{\mbox{lex}}f_s(\alpha )$ si $\alpha\in A_s$.\smallskip

\noindent (c) La suite $(\mbox{Gr}(f_{s^\frown n}))_n$ converge vers $\mbox{Gr}(f_s)$ si $s\in T'_\xi$, et 
$f_\emptyset = \mbox{Id}_Z$.\smallskip

\noindent (d) $\overline{B_p} = B_p\cup B_i$ si $\xi$ est pair, et $\overline{B_i} = B_p\cup B_i$ 
si $\xi$ est impair. De plus, $\bigcup_{s\in T_\xi,~\vert s\vert <m} \mbox{Gr}(f_s)$ est ferm\'e dans 
$Z\times Z$ pour tout entier $m$.\end{defis}

 Il est d\'emontr\'e le th\'eor\`eme suivant dans [Le2] (cf th\'eor\`eme 2.3) :

\begin{thm} Soit $\xi$ un ordinal d\'enombrable.\smallskip

\noindent (1) Si $\xi$ est pair non nul, soient $X$ et $Y$ des espaces polonais, et $A$ un 
bor\'elien $\mbox{pot}(\borathree)$ et $\mbox{pot}(\bormthree)$ de $X\times Y$. Les conditions suivantes 
sont \'equivalentes :\smallskip

\noindent (a) Le bor\'elien $A$ n'est pas potentiellement $D_\xi (\boraone)$.\smallskip

\noindent (b) Il existe une situation g\'en\'erale $(Z,T,(g_s)_{s\in T_\xi})$ et $u : Z\rightarrow X$, 
$v : T\rightarrow Y$ injectives continues telles que 
$\overline{B_p} \cap (u\times v)^{-1}(A) = B_p$.\medskip

\noindent (2) Si $\xi$ est impair, soient $X$ et $Y$ des espaces polonais, et $A$ un bor\'elien 
$\mbox{pot}(\borathree)$ et $\mbox{pot}(\bormthree)$ de $X\times Y$. Les conditions suivantes sont 
\'equivalentes :\smallskip

\noindent (a) Le bor\'elien $A$ n'est pas potentiellement $\check D_\xi (\boraone)$.\smallskip

\noindent (b) Il existe une situation g\'en\'erale $(Z,T,(g_s)_{s\in T_\xi})$ et $u : Z\rightarrow X$, 
$v : T\rightarrow Y$ injectives continues telles que $\overline{B_i} \cap (u\times v)^{-1}(A) = B_i$.\end{thm}

\noindent\bf Notations.\rm~Soit $(\omega^\omega,(f_s)_{s\in T_\xi})$ une bonne situation. On 
d\'efinit une relation $\cal R$ sur $\omega^{<\omega}$ comme en 2 :
$$s~{\cal R}~t~\Leftrightarrow ~\vert s\vert =\vert t\vert ~~\mbox{et}~~(N_s\times N_t)\cap 
\left(\bigcup_{w\in T_\xi} \mbox{Gr}(f_w)\right)\not= \emptyset .$$
$\bullet$ Si $s~{\cal R}~t$, on pose 
$$m(s,t) := \mbox{min}\{m\in\omega~/~\exists~w\in T_\xi~~\vert w\vert =m~~
\mbox{et}~~(N_s\times N_t)\cap \mbox{Gr}(f_w)\not= \emptyset\}.$$

\vfill\eject

\noindent $\bullet$ On pose ensuite 
$$\begin{array}{ll} 
A_s\!\!\!\! 
& := \left\{\alpha\in\omega^\omega~/~\forall~i^{\leq \vert s\vert }_{\geq 1}~~\alpha (N(s\lceil i))=1\right\}\cr 
&\cr 
B_s\!\!\!\! 
& := \left\{\alpha\in\omega^\omega~/~\exists~z\in\omega^{N(s)+1}~~\left\{\!\!\!\!\!\!\!\!
\begin{array}{ll} 
& \forall~i^{\leq \vert s\vert }_{\geq 1}~~z(N(s\lceil i))=1~\mbox{et}~~\forall~p\leq N(s)\cr\cr 
& \alpha (p) = \left\{\!\!\!\!\!\!\!\!
\begin{array}{ll} 
& N(z\lceil (p+1))~\mbox{si}~\exists~i^{\leq \vert s\vert }_{\geq 1}~~p=N(s\lceil i)\mbox{,}\cr &\cr
& z(p)~\mbox{sinon.}
\end{array}
\!\!\!\!\!\!\right.
\end{array}
\right.\right\}\cr\cr
f_s\!\!\!\! 
& :\left\{\!\!
\begin{array}{ll} 
A_s\!\!\!\! 
& \rightarrow B_s \cr\cr 
\alpha 
& \mapsto\left\{\!\!
\begin{array}{ll} 
\omega\!\!\!\! 
& \rightarrow \omega\cr\cr  
p 
& \mapsto \left\{\!\!\!\!\!\!\!\!
\begin{array}{ll} 
& N(\alpha\lceil (p+1))~{\rm si}~\exists~i^{\leq \vert s\vert }_{\geq 1}~~p=N(s\lceil i)\mbox{,}\cr\cr 
& \alpha (p)~\mbox{sinon.} 
\end{array}
\right.
\end{array}
\right.
\end{array}
\right.
\end{array}$$
$\bullet$ Enfin, on pose 
$${\cal A}_n := \left\{\!\!\!\!\!\!\!\!
\begin{array}{ll}
& \{1\}~~\mbox{si}~~n=0~~\mbox{ou}~~n\notin\mbox{Im}(N)\mbox{,}\cr
& \{1\}\cup\{N(s^\frown 1)~/~s\in \Pi_{i<n}~{\cal A}_i~\mbox{et}~
\forall~0<i<\vert N^{-1}(n)\vert ~~s(N(N^{-1}(n)\lceil i)) = 1\}~\mbox{sinon,}
\end{array}\right.$$
$$Z_0 := \Pi_{n\in\omega}~{\cal A}_n.$$
Alors on voit facilement par r\'ecurrence que ${\cal A}_n$ est fini et a au moins deux 
\'el\'ements si $n\!\in\!\mbox{Im}(N)\!\setminus\!\{ 0\}$, de sorte que $Z_0$, muni de la topologie induite par celle de $\omega^\omega$, est hom\'eomorphe \`a $2^\omega$, comme compact m\'etrisable parfait de 
dimension 0 non vide. Il est clair que si $\alpha\in Z_0$ et $\alpha \in A_s$, alors 
$f_s (\alpha )\in Z_0$, de sorte qu'on peut remplacer $\omega^\omega$ par $Z_0$ dans la 
d\'efinition de $f_s$. On note encore $f_s$ cette nouvelle fonction, le contexte 
pr\'ecisant si on travaille dans $\omega^\omega$ ou dans $Z_0$.

\begin{thm} (1) Le couple $(\omega^\omega,(f_s)_{s\in T_\xi})$ est une bonne situation. De plus,\smallskip

\noindent (a) Les classes d'\'equivalence de $\cal E$ sont finies.\smallskip

\noindent (b) Si $c$ est une $\cal T$-cha\^\i ne telle que $\vert c\vert \geq 3$, $c(0)=c(\vert c\vert -1)$, et 
$c(i) \not= c(i+1)$ si $i<\vert c\vert -1$, alors il existe $i<\vert c\vert -2$ tel que 
$c(i) = c(i+2)$.\medskip

\noindent (2) Le couple $(Z_0,(f_s)_{s\in T_\xi})$ est une bonne situation.\end{thm}

\noindent\bf D\'emonstration.\rm\ Les ensembles $A_s$ et $B_s$ sont des 
ouverts-ferm\'es de $\omega^\omega$, clairement. On a que $f_s$ est d\'efinie et bijective, 
et que $B_s \not= \emptyset$ si $A_s \not= \emptyset$. Mais $A_s \not= \emptyset$, 
clairement. Les propri\'et\'es topologiques de $f_s$ sont claires. Si $s$ et $t$ sont deux suites distinctes de $\omega^{<\omega}$, il 
y a deux cas. Ou bien par exemple $s$ d\'ebute $t$ et $\mbox{Gr}(f_s)\cap \mbox{Gr}(f_t) = 
\emptyset$. Ou bien il existe ${i< {\rm min}(\vert s\vert ,\vert t\vert )}$ minimal tel que 
$s(i)\not= t(i)$ ; dans ce cas, $f_s(\alpha)(N[s\lceil (i+1)])$ 
est diff\'erent de $\alpha (N[s\lceil (i+1)])$ alors que 
$f_t(\alpha)(N[s\lceil (i+1)])$ est \'egal \`a 
${\alpha (N[s\lceil (i+1)])}$, par injectivit\'e de $N$. D'o\`u 
$\mbox{Gr}(f_s)\cap \mbox{Gr}(f_t) = \emptyset$ et la condition (b) d'une bonne situation, puisqu'on 
a bien s\^ur $\alpha \leq_{\mbox{lex}} f_s(\alpha )$ si $\alpha\in A_s$. Par ailleurs, la 
condition (a) d'une bonne situation est clairement v\'erifi\'ee. Comme $f_s(\alpha )\in Z_0$ si 
$\alpha\in Z_0\cap A_s$, ces conditions (a) et (b) valent \'egalement dans le cas (2).

\vfill\eject

\noindent $\bullet$ Soit $s\in T'_\xi$, $(\alpha,\beta)\in \omega^\omega\times\omega^\omega$ tels 
que $\forall~n\in\omega~~(\alpha ,\beta )\notin \mbox{Gr}(f_{s^\frown n})$, et 
${(\alpha,\beta) = \displaystyle \lim_{k\rightarrow\infty} {(\alpha_k,\beta_k)}}$, o\`u 
$\beta_k = f_{s^\frown n_k}(\alpha_k )$. Alors on peut supposer que la suite 
d'entiers $(n_k)_k$ tend vers l'infini, pour voir que $\beta = f_s (\alpha )$. 
Comme $\alpha_k\in A_s$, $\alpha\in A_s$. Soit $p\in\omega$ ; alors on trouve un rang $k(p)$ tel que si $k\geq k(p)$, 
$f_{s^\frown n_k}(\alpha_k)(p) = f_s (\alpha_k)(p)$. Ce qui montre que $\beta (p) 
= f_s (\alpha )(p)$. D'o\`u $\overline{\bigcup_{n\in\omega} \mbox{Gr}(f_{s^\frown n})}
\setminus (\bigcup_{n\in\omega} \mbox{Gr}(f_{s^\frown n}))\subseteq \mbox{Gr}(f_s)$.\bigskip

 R\'eciproquement, si $\beta = f_s (\alpha )$, $\forall~n\in\omega~~
(\alpha ,\beta )\notin \mbox{Gr}(f_{s^\frown n})$ et $\alpha\in A_s$. Posons 
$$\alpha_n := \alpha\lceil N(s^\frown n)^\frown 1^\omega .$$ 
On a $\alpha_n\in A_{s^\frown n}$, puisque le fait d'\^etre dans $A_s$ ne d\'epend que des $N(s)+1$ 
premi\`eres coordonn\'ees et que $N(s)<N(s^\frown n)$. On a alors, puisque 
$N(s^\frown n)$ tend vers l'infini avec $n$, que $f_{s^\frown n} (\alpha_n)$ tend 
vers $\beta$. En effet, pour $n$ suffisamment grand, on a ${f_{s^\frown n} 
(\alpha_n)(p) = f_s(\alpha_n) (p) = f_s (\alpha )(p)}$. D'o\`u $(\alpha,\beta)\in 
\overline{\bigcup_{n\in\omega} \mbox{Gr}(f_{s^\frown n})}$, puisque $(\alpha_n)$ tend vers 
$\alpha$. Ceci montre que la suite $(f_s)_{s\in T_\xi}$ v\'erifie la condition (c) d'une 
bonne situation (dans les cas (1) et (2)).\bigskip

\noindent $\bullet$ Montrons que $\bigcup_{s\in T_\xi,~\vert s\vert <m} \mbox{Gr}(f_s)$ est ferm\'e 
dans $\omega^\omega\times\omega^\omega$ si $m$ est entier. Soit donc $(\alpha,\beta)$ dans 
$\omega^\omega\times\omega^\omega$ tel que 
$(\alpha,\beta) = \displaystyle \lim_{k\rightarrow\infty} {(\alpha_k,\beta_k)}$, o\`u 
$\beta_k = f_{s_k} (\alpha_k)$, avec $\vert s_k\vert <m$. Alors on peut supposer que 
$s_k$ est de la forme $s^\frown n_k^\frown u_k$, avec $(n_k)_k$ strictement 
croissante vers l'infini. On montre alors comme pr\'ec\'edemment que $\beta = 
f_s (\alpha )$. D'o\`u le r\'esultat.\bigskip

 Montrons maintenant que $\bigcup_{s\in T_\xi} \mbox{Gr}(f_s)$ est ferm\'e dans 
$\omega^\omega\times\omega^\omega$. Soit 
donc $(\alpha,\beta)$ dans $\omega^\omega\times\omega^\omega$ tel que $(\alpha,\beta) = \displaystyle 
\lim_{k\rightarrow\infty} {(\alpha_k,\beta_k)}$, o\`u $\beta_k = f_{s_k} (\alpha_k)$. 
Par ce qui pr\'ec\`ede, on peut supposer que la suite $(\vert s_k\vert )_k$ cro\^\i t 
strictement vers l'infini. Alors il existe un entier $p$ tel que 
$\{s_k (p)~/~k\in\omega\}$ soit infini. En effet, si tel n'\'etait pas le cas, 
$\{s\in T_\xi~/~\exists~k\in\omega~~s\prec s_k\}$ serait un sous-arbre infini de $T_\xi$, \`a branchements finis, donc il aurait une branche infinie par le lemme de K\" onig. 
Mais ceci contredit la bonne fondation de $T_\xi$. On peut donc supposer qu'il 
existe une suite $s$ telle que $\vert s\vert =p$, $s\prec s_k$ et $(s_k(p))_k$ tende vers 
l'infini. On en d\'eduit comme avant que $\beta = f_s (\alpha )$. D'o\`u le 
r\'esultat. Bien s\^ur, ceci vaut \'egalement dans $Z_0\times Z_0$.\bigskip

\noindent $\bullet$ On a donc $\overline{B_p}\subseteq \bigcup_{s\in T_\xi} \mbox{Gr}(f_s) = 
B_p \cup B_i$. Il reste \`a voir que $B_i\subseteq \overline{B_p}$ pour avoir (d), dans 
le cas o\`u $\xi$ est pair. Soit donc $s\in T_\xi$ de longueur impaire. Comme $s\in 
T'_\xi$, on a la conclusion en utilisant la condition (c). On raisonne de mani\`ere 
analogue si $\xi$ est impair. On a donc montr\'e que les couples nous int\'eressant sont 
des bonnes situations.\bigskip

 Les conditions (a) et (b) de (1) se montrent comme en 2.7, \`a ceci pr\`es que 
${\cal E}(s_0(0)) = \{s_0(0)\}$ vaut dans tous les cas, et que l'injectivit\'e de $N$ 
nous fournit ${c(i_1+1) = c(i_2-1)}$.$\hfill\square$

\begin{thm} Soit $\xi$ un ordinal d\'enombrable pair et 
$(Z,T,(\tilde f_s)_{s\in T_\xi})$ une situation g\'en\'erale. Alors 
il existe $u : \omega^\omega\rightarrow Z$ et $v : \omega^\omega\rightarrow T$ continues telles que  
$\overline{B_p}\cap (u\times v)^{-1}(\tilde B_p) = B_p.$\end{thm} 

\vfill\eject

\noindent\bf D\'emonstration.\rm\ Elle est calqu\'ee sur celle du th\'eor\`eme 2.8. On construit 
des suites $(U_s)_{s\in\omega^{<\omega}}$ et 
$(V_s)_{s\in\omega^{<\omega}}$ d'ouverts, ainsi que $\Phi$ \`a valeurs $T_\xi$ v\'erifiant les conditions 
(i), (ii), et 
$$(iii)~~s~{\cal R}~t~~\Rightarrow~\left\{\!\!\!\!\!\!
\begin{array}{ll} 
& \vert w(s,t)\vert ~\mbox{a~m\^eme~parit\'e~que~et~vaut~au~plus}~m(s,t) \cr 
& V_t = \tilde f_{w(s,t)} [U_s]~\mbox{si}~\tilde A_{w(s,t)}\subseteq Z \cr 
& U_s = \tilde f_{w(s,t)} [V_t]~\mbox{si}~\tilde A_{w(s,t)}\subseteq T
\end{array}
\right.$$
On montre que si $(\alpha ,\beta)$ est dans $B_p$ (resp. $B_i$), alors 
$(u(\alpha ),v(\beta ))$ est dans $\tilde B_p$ (resp. $\tilde B_i$).  Soit 
donc $w\in T_\xi$ tel que  $(\alpha ,\beta)\in \mbox{Gr}(f_w)$ ; alors $\vert w\vert $ est 
paire (resp. impaire), et on peut trouver un entier $m_0$ tel que l'on ait 
${(N_{\alpha\lceil m_0}\times N_{\beta\lceil m_0})\cap 
\bigcup_{s\in T_\xi,~\vert s\vert  < \vert w\vert } \mbox{Gr}(f_s) = \emptyset}$. Par (iii), 
$$m(s,t) = m(\alpha\lceil m_0,\beta\lceil m_0) = \vert w\vert $$ 
et $\vert \Phi(s,t)\vert $ ont 
m\^eme parit\'e. Comme $\vert w\vert $ est paire (resp. impaire), $\vert \Phi(s,t)\vert $ aussi. 
On a alors que $(u(\alpha ),v(\beta ))$ est dans $G(\tilde f_{\Phi(s,t)})\subseteq\tilde B_p$
 (resp. $\tilde B_i$).\bigskip
 
 La suite de la preuve est identique \`a celle du th\'eor\`eme 2.8, le cas 2 non compris, \`a ceci 
pr\`es que le premier alin\'ea de la condition (3) devient 
$$\vert w(k,l)\vert ~\mbox{a~m\^eme~parit\'e~que~et~vaut~
au~plus}~m(z_{\phi^{-1}(k)},z_{\phi^{-1}(l)}).$$
\bf Cas 2.\rm\ $o=p+1$.\bigskip
 
\noindent 2.1. $z_{\phi^{-1}(m)}~{\cal R}~z_{\phi^{-1}(n)}$.\bigskip

 Soit $w\in T_\xi$ tel que 
$(N_{z_{\phi^{-1}(m)\lceil p}}\times 
N_{z_{\phi^{-1}(n)\lceil p}})\cap \mbox{Gr}(f_w)\not=\emptyset$. On peut supposer que 
$$\vert w\vert  = m(z_{\phi^{-1}(m)\lceil p},z_{\phi^{-1}(n)\lceil p}).$$ 
Par (iii), $\vert w(z_{\phi^{-1}(m)\lceil p},z_{\phi^{-1}(n)\lceil p})\vert $ et $\vert w\vert $ ont m\^eme 
parit\'e et $\vert w(z_{\phi^{-1}(m)\lceil p},z_{\phi^{-1}(n)\lceil p})\vert \leq\vert w\vert $. 
Il y a deux cas.\bigskip

\noindent 2.1.1. $\tilde A_{w(z_{\phi^{-1}(m)\lceil p},z_{\phi^{-1}(n)\lceil p})}\subseteq Z$.\bigskip

 Par la condition (3), on a 
$$V_{z_{\phi^{-1}(n)\lceil p}} = \tilde f_{w(z_{\phi^{-1}(m)\lceil p},z_{\phi^{-1}(n)\lceil p})}
[U_{z_{\phi^{-1}(m)\lceil p}}].$$
 On en d\'eduit que $(U_{z_{\phi^{-1}(m)}}^{n-1}\times V_{z_{\phi^{-1}(n)\lceil p}})
\cap G(\tilde f_{w(z_{\phi^{-1}(m)\lceil p},z_{\phi^{-1}(n)\lceil p})})\not=\emptyset$. Comme on a 
$${G(\tilde f_s) = \overline{\bigcup_{n\in\omega} G(\tilde f_{s^\frown n})}\setminus 
\left(\bigcup_{n\in\omega} G(\tilde f_{s^\frown n})\right)}$$ 
si $s\in T'_\xi$, on peut trouver $t\in \omega^{<\omega}$ lexicographiquement minimale telle que 
$$(U_{z_{\phi^{-1}(m)}}^{n-1}\times V_{z_{\phi^{-1}(n)\lceil p}})
\cap G(\tilde f_{w(z_{\phi^{-1}(m)\lceil p},z_{\phi^{-1}(n)\lceil p})^\frown t})
\not=\emptyset\mbox{,}$$ 
et aussi telle que $\vert w(z_{\phi^{-1}(m)\lceil p},z_{\phi^{-1}(n)\lceil p})^\frown t\vert $ soit de m\^eme 
parit\'e que et vaille au plus 
$$m(z_{\phi^{-1}(m)},z_{\phi^{-1}(n)}).$$ 
On pose alors $\Phi(z_{\phi^{-1}(m)},z_{\phi^{-1}(n)}) := 
w(z_{\phi^{-1}(m)\lceil p},z_{\phi^{-1}(n)\lceil p})^\frown t$. Il y a alors 2 cas.

\vfill\eject

\noindent 2.1.1.1. $\tilde A_{w(m,n)} \subseteq Z$.\bigskip

\noindent 2.1.1.1.1. $\tilde A_\emptyset \subseteq Z$.\bigskip

 On a alors $U_{z_{\phi^{-1}(n)\lceil p}} \cap \tilde 
f_\emptyset^{-1}(\tilde f_{w(m,n)}
[U_{z_{\phi^{-1}(m)}}^{n-1}\cap \tilde A_{w(m,n)}])\not=\emptyset$, et on raisonne 
comme en 1.1.1.\bigskip

\noindent 2.1.1.1.2. $\tilde A_\emptyset \subseteq T$.\bigskip

 On a alors $\tilde f_{w(m,n)}[U_{z_{\phi^{-1}(m)}}^{n-1}\cap \tilde A_{w(m,n)}]
\cap \tilde f_\emptyset^{-1}(U_{z_{\phi^{-1}(n)}\lceil p})\not=\emptyset$, et on 
raisonne comme en 1.1.2.\bigskip

\noindent 2.1.1.2. $\tilde A_{w(m,n)} \subseteq T$.\bigskip

\noindent 2.1.1.2.1. $\tilde A_\emptyset \subseteq Z$.\bigskip

 On a alors $U_{z_{\phi^{-1}(n)\lceil p}} \cap \tilde 
f_\emptyset^{-1}(\tilde f_{w(m,n)}^{-1}(
U_{z_{\phi^{-1}(m)}}^{n-1}))\not=\emptyset$, et on raisonne comme en 1.2.1.\bigskip

\noindent 2.1.1.2.2. $\tilde A_\emptyset \subseteq T$.\bigskip

 On a alors $\tilde f_{w(m,n)}^{-1}(U_{z_{\phi^{-1}(m)}}^{n-1})
\cap \tilde f_\emptyset^{-1}(U_{z_{\phi^{-1}(n)}\lceil p})\not=\emptyset$, et on 
raisonne comme en 1.2.2.\bigskip

\noindent 2.1.2. $\tilde A_{w(z_{\phi^{-1}(m)\lceil p},z_{\phi^{-1}(n)\lceil p})}
\subseteq T$.\bigskip

 Par la condition (3), on a $U_{z_{\phi^{-1}(m)\lceil p}} = 
\tilde f_{w(z_{\phi^{-1}(m)\lceil p},z_{\phi^{-1}(n)\lceil p})}
[V_{z_{\phi^{-1}(n)\lceil p}}]$. On conclut alors comme en 2.1.1.\bigskip

\noindent 2.2. $z_{\phi^{-1}(n)}~{\cal R}~z_{\phi^{-1}(m)}$.\bigskip

 On raisonne comme en 2.1 : il y a deux cas.\bigskip
 
\noindent 2.2.1. $\tilde A_{w(z_{\phi^{-1}(n)\lceil p},z_{\phi^{-1}(m)\lceil p})}
\subseteq Z$.\bigskip

 Par la condition (3), on a $V_{z_{\phi^{-1}(m)\lceil p}} = 
\tilde f_{w(z_{\phi^{-1}(n)\lceil p},z_{\phi^{-1}(m)\lceil p})}
[U_{z_{\phi^{-1}(n)\lceil p}}]$. On en d\'eduit que $(U_{z_{\phi^{-1}(n)}\lceil p}\times V_{z_{\phi^{-1}(m)}}^{n-1})
\cap G(\tilde f_{w(z_{\phi^{-1}(n)\lceil p},z_{\phi^{-1}(m)\lceil p})})\not=\emptyset$. 
Comme 
$${G(\tilde f_s) = 
\overline{\bigcup_{n\in\omega} G(\tilde f_{s^\frown n})}\setminus 
\left(\bigcup_{n\in\omega} G(\tilde f_{s^\frown n})\right)}$$ 
si $s\in T'_\xi$, on peut trouver $t'\in \omega^{<\omega}$ lexicographiquement minimale telle que 
$${(U_{z_{\phi^{-1}(n)}\lceil p}\!\times\! V_{z_{\phi^{-1}(m)}}^{n-1})
\!\cap\! G(\tilde f_{w(z_{\phi^{-1}(n)\lceil p},z_{\phi^{-1}(m)\lceil p})^\frown t'})
\!\not=\!\emptyset}\mbox{,}$$ 
et aussi telle que $\vert w(z_{\phi^{-1}(n)\lceil p},z_{\phi^{-1}(m)\lceil p})^\frown t'\vert $ soit de m\^eme 
parit\'e que et vaille au plus 
$$m(z_{\phi^{-1}(n)},z_{\phi^{-1}(m)}).$$ 
On pose alors 
$\Phi(z_{\phi^{-1}(n)},z_{\phi^{-1}(m)}) := 
w(z_{\phi^{-1}(n)\lceil p},z_{\phi^{-1}(m)\lceil p})^\frown t'$. Il y a alors 2 cas.

\vfill\eject

\noindent 2.2.1.1. $\tilde A_{w(n,m)} \subseteq Z$.\bigskip

\noindent 2.2.1.1.1. $\tilde A_\emptyset \subseteq Z$.\bigskip

 On a alors $U_{z_{\phi^{-1}(n)\lceil p}} \cap \tilde f_{w(n,m)}^{-1}(
V_{z_{\phi^{-1}(m)}}^{n-1})\not=\emptyset$, et on raisonne 
comme en 1.3.1.\bigskip

\noindent 2.2.1.1.2. $\tilde A_\emptyset \subseteq T$.\bigskip

 On a alors 
$\tilde f_\emptyset^{-1}(\tilde f_{w(n,m)}^{-1}(V_{z_{\phi^{-1}(m)}}^{n-1}))
\cap V_{z_{\phi^{-1}(n)}\lceil p}\not=\emptyset$, et on 
raisonne comme en 1.3.2.\bigskip

\noindent 2.2.1.2. $\tilde A_{w(n,m)} \subseteq T$.\bigskip

\noindent 2.2.1.2.1. $\tilde A_\emptyset \subseteq Z$.\bigskip

 On a alors 
$\tilde f_{w(n,m)}[V_{z_{\phi^{-1}(m)}}^{n-1}\cap \tilde A_{w(n,m)}]\not=\emptyset$, 
et on raisonne comme en 1.4.1.\bigskip

\noindent 2.2.1.2.2. $\tilde A_\emptyset \subseteq T$.\bigskip

 On a alors $\tilde f_\emptyset^{-1}(\tilde f_{w(n,m)}[V_{z_{\phi^{-1}(m)}}^{n-1}\cap 
\tilde A_{w(n,m)}])\cap V_{z_{\phi^{-1}(n)}\lceil p}\not=\emptyset$, et on raisonne 
comme en 1.4.2.\bigskip

\noindent 2.2.2. $\tilde A_{w(z_{\phi^{-1}(n)\lceil p},z_{\phi^{-1}(m)\lceil p})}
\subseteq T$.\bigskip

 Par la condition (3), on a $U_{z_{\phi^{-1}(n)\lceil p}} = 
\tilde f_{w(z_{\phi^{-1}(n)\lceil p},z_{\phi^{-1}(m)\lceil p})}
[V_{z_{\phi^{-1}(m)\lceil p}}]$. On conclut alors comme en 2.2.1.$\hfill\square$

\begin{thm} Soit $\xi$ un ordinal d\'enombrable.\medskip

\noindent (1) Si $\xi$ est pair, il existe un bor\'elien $A_\xi$ de $2^\omega\times 2^\omega$ tel que 
pour tous espaces polonais $X$ et $Y$, et pour tout bor\'elien $A$ de $X\times Y$ qui est 
$\mbox{pot}(\borathree)$ et $\mbox{pot}(\bormthree)$, on a l'\'equivalence entre les conditions suivantes :\smallskip

\noindent (a) Le bor\'elien $A$ n'est pas $\mbox{pot}(D_\xi (\boraone))$.\smallskip

\noindent (b) Il existe des fonctions continues $u : 2^\omega \rightarrow X$ et 
$v : 2^\omega \rightarrow Y$ telles que $\overline{A_\xi} \cap (u\times v)^{-1}(A) = 
A_\xi$.\medskip

\noindent (2) Si $\xi$ est impair, il existe un bor\'elien $A_\xi$ de $2^\omega\times 2^\omega$ tel que 
pour tous espaces polonais $X$ et $Y$, et pour tout bor\'elien $A$ de $X\times Y$ qui est 
$\mbox{pot}(\borathree)$ et $\mbox{pot}(\bormthree)$, on a l'\'equivalence entre les conditions suivantes :\smallskip

\noindent (a) Le bor\'elien $A$ n'est pas $\mbox{pot}(\check D_\xi (\boraone))$.\smallskip

\noindent (b) Il existe des fonctions continues $u : 2^\omega \rightarrow X$ et 
$v : 2^\omega \rightarrow Y$ telles que $\overline{A_\xi} \cap (u\times v)^{-1}(A) = 
A_\xi$.\end{thm}

\vfill\eject

\noindent\bf D\'emonstration.\rm\ Soit $\Phi : 2^\omega \rightarrow Z_0$ 
un hom\'eomorphisme. On montre le th\'eor\`eme dans le cas o\`u $\xi$ est pair, l'autre cas \'etant 
analogue. On pose 
$$A_\xi := (\Phi\times\Phi )^{-1}(B_p).$$
$\bullet$ Si $\xi=0$, $\overline{B_p}=B_p$. Si $A$ est non vide, soit 
$(x,y)\in A$, et $u$ (resp. $v$) l'application constante identique \`a $x$ 
(resp. $y$). Alors $u$ et $v$ sont continues, et $\overline{A_0} = A_0 \subseteq 
(u\times v)^{-1}(A)$. R\'eciproquement, $A_0\subseteq (u\times v)^{-1}(A)$, donc $A$ est 
non vide et $A\notin D_0(\boraone)$. Dans la suite, on supposera donc $\xi \not= 0$.\bigskip

\noindent $\bullet$ Appliquons le th\'eor\`eme 3.2 \`a $X=Y=2^\omega$, $A=A_\xi$, $Z=T=Z_0$, 
et $u=v=\Phi^{-1}$. Ce th\'eor\`eme s'applique car\bigskip

\noindent - $A_\xi$ a ses coupes d\'enombrables, donc est $\mbox{pot}(\borathree)$ et $\mbox{pot}(\bormthree)$.\bigskip

\noindent - $(Z_0,(f_s)_{s\in T_\xi})$ est une bonne situation (par 3.3), donc 
$(Z_0,Z_0,(f_s)_{s\in T_\xi})$ est une situation g\'en\'erale.\bigskip

\noindent - $B_p = (u\times v)^{-1}(A_\xi) = (u\times v)^{-1}(A_\xi)\cap\overline{B_p}$.\bigskip 

 On a alors que $A_\xi$ n'est pas $\mbox{pot}(D_\xi(\boraone))$, ce qui montre que (b) 
implique (a) (la parit\'e de $\xi$ fait que la classe $D_\xi(\boraone)$ est stable par 
intersection avec les ferm\'es).\bigskip

\noindent $\bullet$ R\'eciproquement, si $A$ n'est pas potentiellement $D_\xi(\boraone)$, 
le th\'eor\`eme 3.2 nous fournit une situation g\'en\'erale $(Z,T,(\tilde f_s)_{s\in T_\xi})$ et des injections 
continues $\tilde u : Z\rightarrow X$ et $\tilde v : T\rightarrow Y$ telles que 
$$\overline{\tilde B_p} \cap (\tilde u\times\tilde v)^{-1}(A) = \tilde B_p.$$ 
Le th\'eor\`eme 3.4 nous fournit des fonctions continues 
$u':\omega^\omega\rightarrow Z$ et $v':\omega^\omega\rightarrow T$ telles que 
$$\overline{B_p}\cap (u'\times v')^{-1}(\tilde B_p) = B_p.$$ 
Soit $\Psi: Z_0\rightarrow \omega^\omega$ l'injection canonique. On pose 
$$u:=\tilde u\circ~u'\circ~\Psi\circ~\Phi\mbox{,}~~v:=\tilde v\circ~v'\circ~\Psi\circ~\Phi .$$
Alors $u$ et $v$ sont clairement continues, et on a clairement 
$A_\xi\subseteq (u\times v)^{-1}(A)$. Si $(\alpha,\beta)\in \overline{A_\xi}\setminus 
A_\xi$, $(\Phi(\alpha ),\Phi(\beta ))\in \overline{B_p}\setminus B_p = B_i$, donc 
$(u'\circ~\Psi\circ~\Phi(\alpha ), v'\circ~\Psi\circ~\Phi(\beta ))\in \tilde B_i
\subseteq \overline{\tilde B_p}$, et $(u(\alpha ),v(\beta ))$ n'est pas dans $A$.$\hfill\square$\bigskip

\bf\noindent Remarques.\rm~(a) L'\'enonc\'e de ce th\'eor\`eme fournit un exemple de test $A_\xi$ 
dans $2^\omega \times 2^\omega$ ; mais la preuve nous donne aussi un test dans 
$Z_0\times Z_0$, et un autre dans $\omega^\omega \times \omega^\omega$.\bigskip

\noindent (b) Le th\'eor\`eme 3.5 fournit une caract\'erisation des ensembles non potentiellement 
$D_\xi(\boraone)$ pour $\xi$ pair, et non potentiellement $\check D_\xi(\boraone)$ pour 
$\xi$ impair. Mais un simple passage au compl\'ementaire nous fournit une caract\'erisation des 
ensembles non potentiellement $\check D_\xi(\boraone)$ pour $\xi$ pair, et non 
potentiellement $D_\xi(\boraone)$ pour $\xi$ impair :

\vfill\eject

\begin{thm} Soit $\xi$ un ordinal d\'enombrable.\medskip

\noindent (1) Si $\xi$ est pair, il existe un bor\'elien $A_\xi$ de $2^\omega\times 2^\omega$ tel que 
pour tous espaces polonais $X$ et $Y$, et pour tout bor\'elien $A$ de $X\times Y$ qui est 
$\mbox{pot}(\borathree)$ et $\mbox{pot}(\bormthree)$, on a l'\'equivalence entre les conditions suivantes :\smallskip

\noindent (a) Le bor\'elien $A$ n'est pas $\mbox{pot}(\check D_\xi (\boraone))$.\smallskip

\noindent (b) Il existe $u : 2^\omega \rightarrow X$ et $v : 2^\omega \rightarrow Y$ continues telles que 
$\overline{A_\xi} \cap (u\times v)^{-1}(A) = \overline{A_\xi}\setminus A_\xi$.\medskip

\noindent (2) Si $\xi$ est impair, il existe un bor\'elien $A_\xi$ de $2^\omega\times 2^\omega$ tel que 
pour tous espaces polonais $X$ et $Y$, et pour tout bor\'elien $A$ de $X\times Y$ qui est 
$\mbox{pot}(\borathree)$ et $\mbox{pot}(\bormthree)$, on a l'\'equivalence entre les conditions suivantes :\smallskip

\noindent (a) Le bor\'elien $A$ n'est pas $\mbox{pot}(D_\xi (\boraone))$.\smallskip

\noindent (b) Il existe $u : 2^\omega \rightarrow X$ et $v : 2^\omega \rightarrow Y$ continues telles que 
$\overline{A_\xi} \cap (u\times v)^{-1}(A) = \overline{A_\xi}\setminus A_\xi$.\end{thm}

En fait, on a $\overline{A_\xi}\setminus A_\xi = (\Phi\times \Phi)^{-1}(B_i)$ si $\xi$ est pair ; si $\xi$ est impair, on a $A_\xi = (\Phi\times \Phi)^{-1}(B_i)$ et $\overline{A_\xi}\setminus A_\xi = (\Phi\times \Phi)^{-1}(B_p)$.

\section{$\!\!\!\!\!\!$ R\'ef\'erences.}

\noindent [Ku]\ \  K. Kuratowski,~\it Topology,~\rm Vol. 1, Academic Press, 1966

\noindent [Le1]\ \ D. Lecomte,~\it Classes de Wadge potentielles et th\'eor\`emes d'uniformisation 
partielle,~\rm Fund. Math.~143 (1993), 231-258

\noindent [Le2]\ \ D. Lecomte,~\it Uniformisations partielles et crit\`eres \`a la Hurewicz dans le plan,~\rm Trans. A.M.S. ~347, 11 (1995), 4433-4460

\noindent [Lo1]\ \ A. Louveau,~\it A separation theorem for $\Ana$ sets,~\rm Trans. A. M. S.~260 (1980), 363-378

\noindent [Lo2]\ \ A. Louveau,~\it Some results in the Wadge hierarchy of Borel sets,~\rm Cabal Sem. 79-81 (A. S. Kechris, D. A. Mauldin, Y. N. Moschovakis, eds), Lect. Notes in Math. 1019 Springer-Verlag (1983), 28-55

\noindent [Lo3]\ \ A. Louveau,~\it Livre \`a para\^\i tre\rm

\noindent [Lo-SR]\ \ A. Louveau and J. Saint Raymond,~\it Borel classes and closed games : Wadge-type and Hurewicz-type results,~\rm Trans. A. M. S.~304 (1987), 431-467

\noindent [Mo]\ \ Y. N. Moschovakis,~\it Descriptive set theory,~\rm North-Holland, 1980

\noindent [SR]\ \ J. Saint Raymond,~\it La structure bor\'elienne d'Effros est-elle standard ?,~\rm 
Fund. Math.~100 (1978), 201-210

\end{document}